\documentclass[10pt,twoside,english]{article}

\newtheorem{Theoreme}{Theorem}

\newtheorem{Proposition}[Theoreme]{Proposition}

\newtheorem{Remarque}[Theoreme]{Remark}

\usepackage{amssymb,amsmath,babel,geometry,latexsym,graphics,tabularx,shapepar,enumerate, pdfsync}
\ifx\optionkeymacros\undefined\else \fi

\catcode`\Œ=\active\defŒ{{\aa}}       
\catcode`\º=\active\defº{\int}        
\catcode`\=\active\def{\c c}        
\catcode`\¶=\active\def¶{\partial}    
\catcode`\Ä=\active\defÄ{\oint}       
\catcode`\Æ=\active\defÆ{\triangle}   
\catcode`\Â=\active\defÂ{\neg}        
\catcode`\µ=\active\defµ{\mu}         
\catcode`\¿=\active\def¿{{\o}}        
\catcode`\¹=\active\def¹{\pi}         
\catcode`\Ï=\active\defÏ{{\oe}}       
\catcode`\§=\active\def§{{\ss}}       
\catcode`\ =\active\def {\dagger}     
\catcode`\Ã=\active\defÃ{\sqrt}       
\catcode`\·=\active\def·{\Sigma}      
\catcode`\Å=\active\defÅ{\approx}     
\catcode`\½=\active\def½{\Omega}      
\catcode`\£=\active\def£{{\it\$}}     
\catcode`\°=\active\def°{\infty}      
\catcode`\¤=\active\def¤{{\S}}        
\catcode`\¦=\active\def¦{{\P}}        
\catcode`\¥=\active\def¥{\bullet}     
\catcode`\»=\active\def»{\leavevmode\raise.585ex\hbox{\b a}}      
\catcode`\¼=\active\def¼{\leavevmode\raise.6ex\hbox{\b o}}        
\catcode`\­=\active\def­{\not=}       
\catcode`\²=\active\def²{\leq}        
\catcode`\³=\active\def³{\geq}        
\catcode`\Ö=\active\defÖ{\div}        
\catcode`\É=\active\defÉ{{\dots}}     
\catcode`\¾=\active\def¾{{\ae}}       
\catcode`\Ç=\active\defÇ{\ll}         
\catcode`\Ò=\active\defÒ{``}          
\catcode`\Á=\active\defÁ{!`}          
\catcode`\¢=\active\def¢{\rlap/c}     
\catcode`\Ô=\active\defÔ{`}           
\catcode`\Õ=\active\defÕ{'}           

\catcode`\=\active\def{{\AA}}       
\catcode`\'=\active\def'{\c C}        
\catcode`\¯=\active\def¯{{\O}}        
\catcode`\¸=\active\def¸{\Pi}         
\catcode`\Î=\active\defÎ{{\OE}}       
\catcode`\®=\active\def®{{\AE}}       
\catcode`\×=\active\def×{\diamond}    
\catcode`\¡=\active\def¡{\accent'27}  
\catcode`\Ó=\active\defÓ{''}          
\catcode`\±=\active\def±{\pm}         
\catcode`\È=\active\defÈ{\gg}         
\catcode`\À=\active\defÀ{?`}          
\catcode`\Ð=\active\defÐ{--}          
\catcode`\Ñ=\active\defÑ{---}         

\catcode`\Š=\active\defŠ{\"a}        
\catcode`\'=\active\def'{\"e}        
\catcode`\•=\active\def•{\"{\i}}     
\catcode`\š=\active\defš{\"o}        
\catcode`\Ÿ=\active\defŸ{\"u}        
\catcode`\Ø=\active\defØ{\"y}        
\catcode`\€=\active\def€{\"A}        
\catcode`\…=\active\def…{\"O}        
\catcode`\†=\active\def†{\"U}        
\catcode`\‡=\active\def‡{\'a}        
\catcode`\Ž=\active\defŽ{\'e}        
\catcode`\'=\active\def'{\'{\i}}     
\catcode`\—=\active\def—{\'o}        
\catcode`\œ=\active\defœ{\'u}        
\catcode`\ƒ=\active\defƒ{\'E}        
\catcode`\ˆ=\active\defˆ{\`a}        
\catcode`\=\active\def{\`e}        
\catcode`\"=\active\def"{\`{\i}}     
\catcode`\˜=\active\def˜{\`o}        
\catcode`\=\active\def{\`u}        
\catcode`\Ë=\active\defË{\`A}        
\catcode`\‹=\active\def‹{\~a}        
\catcode`\–=\active\def–{\~n}        
\catcode`\›=\active\def›{\~o}        
\catcode`\Ì=\active\defÌ{\~A}        
\catcode`\"=\active\def"{\~N}        
\catcode`\Í=\active\defÍ{\~O}        
\catcode`\‰=\active\def‰{\^a}        
\catcode`\=\active\def{\^e}        
\catcode`\"=\active\def"{\^{\i}}     
\catcode`\™=\active\def™{\^o}        
\catcode`\ž=\active\defž{\^u}        

\let\optionkeymacros\null

\usepackage{pdfsync}
\usepackage{hyperref}
\usepackage[all,2cell]{xy} \UseAllTwocells \SilentMatrices
\def\carlitz{{{\tiny \text{Car}}}}

\newcommand{\ZZ}{\mathbb{Z}}
\newcommand{\FF}{\mathbb{F}}
\newcommand{\CC}{\mathbb{C}}
\newcommand{\QQ}{\mathbb{Q}}
\newcommand{\RR}{\mathbb{R}}

\newcommand{\GG}{\mathbb{G}}
\newcommand{\KK}{\mathbb{K}}

\newcommand{\NN}{\mathbb{N}}
\newcommand{\PP}{\mathbb{P}}

\newcommand{\und}{\underline}
\newcommand{\ord}{\text{ord}}

\newcommand\CVD{{\hfill\hfil{\lower 2 pt\hbox{\vrule\vbox to 7pt 
{\hrule width 6pt\vfill\hrule}\vrule}}}\vskip 0.5cm}

\title{An introduction to Mahler's method for transcendence and algebraic independence.}
\author{Federico Pellarin}

\begin{document}

\maketitle
\setcounter{tocdepth}{1}
\tableofcontents

\section{Introduction}

In his mathematical production, Kurt Mahler (1903-1988) introduced
entire new subjects. One of them, perhaps chronologically the first one, 
aimed three fundamental papers \cite{Ma1, Ma2, Ma3} and was, later in 1977, baptised
``Mahler's method" by Loxton and van der Poorten. By this locution we mean a
general method to prove transcendence and algebraic independence 
of values of a certain class of functions by means of the following classical scheme of demonstration whose terminology will be explained in the present text
(\footnote{We borrowed this presentation from Masser's article \cite[p. 5]{Mass}, whose point of view influenced our point of view.}):
\begin{description}
\item (AP) - Construction of {\em auxiliary polynomials},
\item (UP) - Obtaining an {\em upper bound}, by means of {\em analytic estimates},
\item (NV) - Proving the {\em non-vanishing}, by means of {\em zeros estimates},
\item (LB) - Obtaining a {\em lower bound}, by means of {\em arithmetic estimates}.
\end{description}

For example, with Mahler's method, and with the help of the basic theory of {\em heights}, it is possible to show the transcendence of values 
at algebraic complex numbers of transcendental analytic solutions $f(x)\in \mathcal{L}[[x]]$
(with $ \mathcal{L}$ a number field embedded in $\CC$) of functional equations such as
\begin{equation}\label{mahlerfunction}
f(x^d)=R(x,f(x)),\quad d>1,\quad  R\in  \mathcal{L}(X,Y),\end{equation} with $d$ integer, see \cite{Ma1}.

In this text, we will also be interested in analogues of these functions over 
complete, algebraically closed fields other than $\CC$ and for this purpose it will be advantageous
to choose right away an appropriate terminology.
Indeed, in the typical situation we will analyse, there will be a {\em base field}
$\mathcal{K}$, together with a distinguished absolute value that will be denoted by $|\cdot|$.
Over $\mathcal{K}$ there will be other absolute values as well, and a {\em product formula}
will hold. We will consider the completion of $\mathcal{K}$ with respect to $|\cdot|$, its
algebraic closure that will be embedded in its completion $\KK$ with respect to an extension of $|\cdot|$.
The algebraic closure of $\mathcal{K}^{\text{alg.}}$, embedded in $\KK$, will also be 
endowed with a {\em absolute logarithmic height} that will be used to prove transcendence
results. Here, an element of $\KK$ is {\em transcendental} if it does not belong to $\mathcal{K}^{\text{alg.}}$. 

If $ \mathcal{L}$ is a finite extension of $\mathcal{K}$ in $\KK$ and $f\in \mathcal{L}[[x]]$ is 
a formal series solution of (\ref{mahlerfunction}), we will say that $f$ is a {\em Mahler's function 
over $\mathcal{K}$}. If $f$ converges at $\alpha\in\mathcal{K}^{\text{alg.}}\setminus\{0\}$ (for the distinguished absolute value), we
will say that $f(\alpha)\in\KK$ is a {\em Mahler's value} and $\alpha$ is a {\em base point} for this value.
In spite of the generality of this terminology, in this text we will restrict our attention to the base fields
$\QQ,K=\FF_q(\theta)$ and $C(t)$ where $C$ is the completion of an algebraic closure of the completion of $K$ 
for the unique extension of the absolute value defined by $|a|=q^{\deg_\theta a}$, with $a\in K$.

The interest of the method introduced by Mahler in \cite{Ma1} is that it can be generalised, as it was remarked by Mahler himself in \cite{Ma2}, to explicitly produce finitely generated subfields of $\CC$ 
of arbitrarily large transcendence degree.
This partly explains, after that the theory was long-neglected for about forty years,
a regain of interest in it, starting from the late seventies, especially due to the intensive work of Loxton and van der Poorten,
Masser, Nishioka as well as other authors we do not mention here but that are quoted, for example, 
in Nishioka's book \cite{Ni}.

In some sense, Mahler's functions and Siegel's $E$-functions share similar properties; 
large transcendence degree subfields of $\CC$ can also be explicitly constructed 
by the so-called Siegel-Shidlowski theorem on values of Siegel $E$-functions at algebraic
numbers (see Lang's account on the theory in \cite{La}). However, this method makes fundamental use of the fact that $E$-functions
are entire, with finite analytic growth order. This strong assumption 
is not at all required when it is possible to apply Mahler's method,
where the functions involved have natural boundaries for analytic continuation; this is certainly
an advantage that this theory has.
Unfortunately, no complex ``classical constant" (period, special value of exponential function at algebraic 
numbers\ldots)
seems to occur as a complex Mahler's value, as far as we can see. 

More recently, a variety of results by 
Becker, Denis \cite{Becker0,Denis2, Denis1, denis,denis3} and other authors 
changed the aspect of the theory, especially that of Mahler's functions over fields of positive characteristic. It was a
fundamental discovery of Denis, that every period of Carlitz's exponential function
is a Mahler's value, hence providing a new proof of its transcendency. This motivates our choice of terminology; we hope the reader
will not find it too heavy. At least, it will be useful to compare the theory over $\QQ$ and that over $K$.


The aim of this paper is to provide an overview of the theory from its beginning
(transcendence) to its recent development in algebraic independence and its important
excursions in positive characteristic, where it is in ``competition" with 
more recent, and completely different techniques inspired by the theory of Anderson's {\em $t$-motives}
(see, for example, the work of Anderson, Brownawell, Papanikolas, Chieh-Yu Chang, Jing Yu, other authors
\cite{ABP, ChangYu, Papanikolas1} and the related 
bibliographies).

The presentation of the paper essentially follows, in an expanded form, the instructional talk the author gave in the conference
``$t$-motives: Hodge structures, transcendence and other motivic aspects", held
in Banff, Alberta, (September 27 - October 2, 2009). The author is thankful to the organisers of this excellent conference for giving the  opportunity to present 
these topics, and thankful to the Banff Centre for the exceptional environment of working it provided.
The author also wishes to express his gratitude to B. Adamczewski and P. Philippon for 
discussions and hints that helped to improve the presentation of this text, and to P. Bundschuh, H. Kaneko and T. Tanaka
for a description of the algebraic relations involving the functions $L_r$ of Section \ref{examplecomplex} they provided.

Here is what the paper contains. In Section \ref{transcendence}, we give an account of transcendence theory of Mahler's values with $\QQ$ as a base field; this is part of the classical theory, essentially contained in one of the first results by Mahler.
In Section \ref{anotherexample}, we will outline the transcendence theory with, as a base field,
a function field of positive characteristic (topic which is closer to the themes of the conference). Here,
the main two features are some applications to the arithmetic of periods of Anderson's $t$-motives and
some generalisations of results of the literature (cf. Theorem \ref{theo1}).
In Section \ref{algindep}, we first make an overview of known results of algebraic independence over $\QQ$ of Mahler's values, then we describe more recent results in positive characteristic (with the base field $K=\FF_q(\theta)$) and finally, we mention some quantitative aspects. The main features of this section are elementary proofs of two results: one by Papanikolas \cite{Papanikolas1},
describing algebraic dependence relations between certain special values of Carlitz's logarithms,
and another one, by Chieh-Yu Chang and Jing Yu, describing all the algebraic dependence relations of values of Carlitz-Goss 
zeta function at positive integers.

This paper does not contain a complete survey on Mahler's method. For example, 
Mahler's method was also successful in handling several variable
functions. To keep the size of this survey reasonable, we made the arguable decision of not 
describing this part of the theory, concentrating on the theory in one variable, 
which seemed closer to the other themes of the conference. 
\section{Transcendence theory over the base field $\QQ$\label{transcendence}}

\subsection{An example to begin with\label{firstexample}}

The example that follows gives an idea of the method.
We consider the formal series:
$$f_{\text{TM}}(x)=\prod_{n=0}^\infty(1-x^{2^n})=\sum_{n=0}^\infty(-1)^{a_n}x^n\in\ZZ[[x]]$$
$(a_n)_{n\geq0}$ being the Thue-Morse sequence ($a_n$ is the reduction modulo $2$ of the sum of 
the digits of the binary expansion of $n$ and, needless to say, the subscript TM in $f_{\text{TM}}$ stands for ``Thue-Morse").
The formal series $f_{\text{TM}}$ converges in the open unit ball $B(0,1)$ to an analytic function and satisfies the functional equation
\begin{equation}\label{functionaleqf}
f_{\text{TM}}(x^2)=\frac{f_{\text{TM}}(x)}{1-x}, 
\end{equation}
(in (\ref{mahlerfunction}), $d=2$ and $R=\frac{Y}{1-X}$ so that $f_{\text{TM}}$ is a Mahler's function. In all the following,
we fix an embedding in $\CC$ of the algebraic closure $\QQ^{\text{alg.}}$ of $\QQ$.
We want to prove:

\begin{Theoreme} For all $\alpha\in \QQ^{\text{alg.}}$ with $0<|\alpha|<1$, $f_{\text{TM}}(\alpha)$ is transcendental.\label{firstexamplebis}\end{Theoreme}

This is a very particular case of a result of Mahler \cite{Ma1} reproduced as Theorem \ref{mahler_nishioka}
in the present paper.
The proof, contained in \ref{proofofexample} uses properties of {\em Weil's logarithmic absolute height} 
reviewed in \ref{height}. It will also use the property that $f_{\text{TM}}$ is transcendental over $\CC(x)$,
proved below in \ref{ftranscendental}.

\subsubsection{Transcendence of $f_{\text{TM}}$.\label{ftranscendental}}
The transcendence of $f_{\text{TM}}$ over $\CC(x)$ can be checked in several ways.
A first way to proceed appeals to P\'olya-Carlson Theorem (1921), (statement and proof can be found on p. 265 of \cite{Re}). It
says that a given formal series $\phi\in\ZZ[[x]]$ converging with radius of convergence $1$, either has $\{z\in\CC,|z|=1\}$ as
natural boundary for holomorphy, or can be extended to a rational function of the form
$\frac{P(x)}{(1-x^m)^n}$, with $P\in\ZZ[x]$. To show that $f_{\text{TM}}$ is transcendental, it suffices to prove that $f_{\text{TM}}$ is not of the 
form above, which is evident from the functional equation (\ref{functionaleqf}), which implies that $f_{\text{TM}}$ has bounded integral coefficients.
Indeed, if rational, $f_{\text{TM}}$ should have
ultimately periodic sequence of the coefficients. However, it is well known 
(and easy to prove, see \cite[Chapter 5, Proposition 5.1.2]{FoBe})
that this is not the case for the Thue-Morse sequence.

Another way to check the transcendence of $f_{\text{TM}}$ is that suggested in Nishioka's paper \cite{Ni3}.
Assuming that $f_{\text{TM}}$ is algebraic, the field $F=\CC(x,f_{\text{TM}}(x))$ is an algebraic extension of $\CC(x)$ of degree, say $n$, and we want to prove that this degree is $1$.
It is possible to contradict this property observing that the extension $F$ of $\CC(x^d)$ ramifies at the places $0$ and $\infty$ only and applying Riemann-Hurwitz formula. Hence, $f_{\text{TM}}$ is rational and
we know already from the lines above how to exclude this case.


\subsubsection{A short account on heights.\label{height}}

Here we closely follow Lang \cite[Chapter 3]{La} and Waldschmidt \cite[Chapter 3]{Wa}.
Let $L$ be a number field. The {\em absolute logarithmic height}
$h(\alpha_0:\cdots:\alpha_n)$
of a projective point $(\alpha_0:\cdots:\alpha_n)\in\PP_n(L)$ is the following weighted average of 
logarithms of absolute values:
$$h(\alpha_0:\cdots:\alpha_n)=\frac{1}{[L:\QQ]}\sum_{v\in M_L}d_v\log\max\{|\alpha_0|_v,\ldots,|\alpha_n|_v\},$$
where $v$ runs over a complete set $M_L$ of non-equivalent places of $L$, where $d_v=[L_v:\QQ_p]$ with $v|_\QQ=p$ the {\em local degree} at the place $v$
(one then writes that $v|p$)
($L_v,\QQ_p$ are completions of $L,\QQ$ at the respective places so that if $v|\infty$,
$L_v=\RR$ or $L_v=\CC$ according to whether the place $v$ is real or complex), and where $|\cdot|_v$ denotes, for all $v$, an element of $v$ chosen in such a way that the following {\em product formula} holds:
\begin{equation}
\prod_{v\in M_L}|\alpha|_v^{d_v}=1,\quad \alpha\in L^\times,
\label{productformula}\end{equation}
where we notice that only finitely many factors of this product are distinct from $1$.
A common way to normalise the $|\cdot|_v$'s is to set $|x|_v=x$ if $x\in\QQ$, $x>0$, and $v|\infty$,
and $|p|_v=1/p$ if $v|p$. 

This formula implies that $h$ does not depend on the choice of the number field $L$,
so that we have a well defined function $$h:\PP_n(\QQ^{\text{alg.}})\rightarrow\RR_{\geq0}.$$

If $n=1$ we also write $h(\alpha)=h(1:\alpha)$. For example, we have $h(p/q)=h((1:p/q))=\log\max\{|p|,|q|\}$
if $p,q$ are relatively prime and $q\neq 0$. With the convention $h(0):=0$, this defines a function 
\begin{equation}\label{heightonevariable}h:\QQ^{\text{alg.}}\rightarrow\RR_{\geq0}
\end{equation} satisfying, for $\alpha,\beta\in\QQ^{\text{alg.}}{}^\times$:
\begin{eqnarray*}
h(\alpha+\beta)&\leq &h(\alpha)+h(\beta)+\log 2,\\
h(\alpha\beta)&\leq&h(\alpha)+h(\beta),\\
h(\alpha^n)&=&|n|h(\alpha),\quad n\in\ZZ.
\end{eqnarray*}
More generally, if $P\in\ZZ[X_1,\ldots,X_n]$ and if $\alpha_1,\ldots,\alpha_n$ are in $\QQ^{\text{alg.}}$,
\begin{equation}\label{eqlio}h(P(\alpha_1,\ldots,\alpha_n))\leq\log L(P)+\sum_{i=1}^n(\deg_{X_i}P)h(\alpha_i),
\end{equation}
where $L(P)$ denotes the {\em length} of $P$, that is,
the sum of the absolute values of the coefficients of $P$.
Proofs of these properties are easy collecting metric information at every place. More details can be found in \cite[Chapter 3]{Wa}.

\medskip

{\em Liouville's inequality}, a sort of ``fundamental theorem of transcendence", reads as follows. 
Let $L$ be a number field, $v$ an archimedean place of $L$, $n$ an integer. For $i=1,\ldots,n$, let $\alpha_i$ be an element of $L$.
Further, let $P$ be a polynomial in $n$ variables $X_1,\ldots,X_n$, with coefficients in $\ZZ$, which does not vanish at the point $(\alpha_1,\ldots,\alpha_n)$.
Assume that $P$ is of degree at most $N_i$ with respect to the variable $X_i$. Then,
$$
\log|P(\alpha_1,\ldots,\alpha_n)|_v\geq -([L:\QQ]-1)\log L(P)-[L:\QQ]\sum_{i=1}^nN_ih(\alpha_i).
$$
The proof of this inequality is again a simple application of product formula (\ref{productformula}): \cite[Section 3.5]{Wa}. It implies that for $\beta\in\QQ^{\text{alg.}}$, $\beta\neq0$,
\begin{equation}\label{liouville}\log|\beta|\geq-[L:\QQ]h(\beta).\end{equation} This inequality suffices for most of 
the arithmetic purposes of this paper
(again, see \cite{Wa} for the details of these basic tools).

\subsubsection{Proof of Theorem \ref{firstexamplebis}\label{proofofexample}}

\noindent\emph{Step (AP)}. For all $N\geq 0$, we choose a polynomial $P_N\in\QQ[X,Y]\setminus\{0\}$ of degree $\leq N$ in both $X$ and $Y$, such that the {\em order  of vanishing} $\nu(N)$ at $x=0$ of the formal series
$$F_N(x):=P_N(x,f_{\text{TM}}(x))=c_{\nu(N)}x^{\nu(N)}+\cdots\quad (c_{\nu(N)}\neq0).$$ (not identically zero because $f_{\text{TM}}$ is transcendental by \ref{ftranscendental}), is $\geq N^2$.

The existence of $P_N$ follows from the existence of a non-trivial 
solution of a homogeneous linear system with $N^2$ linear equations defined over $\QQ$ in $(N+1)^2$
indeterminates. We will not need to control the size of the coefficients of $P_N$ and this is quite unusual in transcendence theory.

\medskip

\noindent\emph{Step (NV)}. Let $\alpha$ be an algebraic number such that $0<|\alpha|<1$ and let us suppose by contradiction that $f_{\text{TM}}(\alpha)$ is also algebraic, so that there exists 
a number field $L$ containing at once $\alpha$ and $f_{\text{TM}}(\alpha)$. Then, by the functional equation (\ref{functionaleqf}), for all $n\geq 0$,
$$F_N(\alpha^{2^{n+1}})=P_N(\alpha^{2^{n+1}},f_{\text{TM}}(\alpha^{2^{n+1}}))=P_N\left(\alpha^{2^{n+1}},\frac{f_{\text{TM}}(\alpha)}{(1-\alpha)\cdots(1-\alpha^{2^n})}\right)\in L.$$
We know that $F_N(\alpha^{2^{n+1}})\neq0$ for all $n$ big enough depending on $N$ and $\alpha$; indeed,
$f_N$ is not identically zero and analytic at $0$.

\medskip

\noindent\emph{Step (UB)}. 
Writing the expansion of $f_N$ at $0$ $$f_N(x)=\sum_{m\geq \nu(N)}c_mx^m=x^{\nu(N)}\left(c_{\nu(N)}+\sum_{i\geq 1}c_{\nu(N)+i}x^i\right)$$
with the leading coefficient $c_{\nu(N)}$ which is a non-zero rational integer (whose size we do not control), we
see that for all $\epsilon>0$, if $n$ is big enough depending on $N,\alpha$ and $\epsilon$:
$$\log|F_N(\alpha^{2^{n+1}})|\leq \log|c_{\nu(N)}|+2^{n+1}\nu(N)\log|\alpha|+\epsilon.$$

\medskip

\noindent\emph{Step (LB)}. 
At once, by (\ref{eqlio}) and (\ref{liouville}),
\begin{eqnarray*}
\lefteqn{\log|F_N(\alpha^{2^{n+1}})|\geq}\\ & \geq &-[L:\QQ](L(P_N)+Nh(\alpha^{2^{n+1}})+Nh(f(\alpha)/(1-\alpha)\cdots(1-\alpha^{2^n})))\\
&\geq &-[L:\QQ](L(P_N)+N2^{n+1}h(\alpha)+Nh(f_{\text{TM}}(\alpha))+\sum_{i=0}^nh(1-\alpha^{2^n}))\\
&\geq &-[L:\QQ](L(P_N)+N2^{n+2}h(\alpha)+Nh(f_{\text{TM}}(\alpha))+(n+1)\log2).
\end{eqnarray*}

\medskip

The four steps allow to conclude:
for all $n$ big enough,
\begin{eqnarray*}
\lefteqn{2^{-n-1}\log|c_{\nu(N)}|+\nu(N)\log|\alpha|+2^{-n-1}\epsilon\geq}\\
&\geq &-[L:\QQ](L(P_N)2^{-n-1}+2Nh(\alpha)+N2^{-n-1}h(f_{\text{TM}}(\alpha))+(n+1)2^{-n-1}\log2).
\end{eqnarray*}
Letting $n$ tend to infinity and using that $\nu(N)\geq N^2$ (recall that $\log|\alpha|$ is negative),
we find the inequality
$$N\log|\alpha|\geq-2[L:\QQ]h(\alpha).$$
But the choice of the ``auxiliary" polynomial $P_N$ can be done for {\em every} $N>0$. 
With
\begin{equation}\label{boundN}
N>\frac{2[L:\QQ]h(\alpha)}{|\log|\alpha||},
\end{equation}
we encounter a contradiction.

\subsubsection{A more general result\label{moregeneralres}}

For $R=N/D\in\CC(X,Y)$ with $N,D$ relatively prime polynomials in $\CC(X)[Y]$, we write $h_Y(R):=\max\{\deg_YN,\deg_YD\}$. With the arguments above, the reader can be easily prove the following theorem originally due to 
Mahler \cite{Ma1}.

\begin{Theoreme}[Mahler]\label{mahler_nishioka}
Let $L\subset\CC$ be a number field, $R$ be an element of $L(X,Y)$,
$d>1$ an integer such that $h_Y(R)<d$. Let $f\in L[[x]]$ be a transcendental formal series such that, in $L((x))$,
$$f(x^d)=R(x,f(x)).$$
Let us suppose that $f$ converges for $x\in\CC$ with $|x|<1$. Let $\alpha$ be an element 
of $L$ such that $0<|\alpha|<1$.

Then, for all $n$ big enough, $f(\alpha^{d^n})$ is transcendental over $\QQ$. 
\end{Theoreme}

Obviously, for all $n$ big enough, $L(f(\alpha^{d^n}))^{\text{alg.}}=L(f(\alpha^{d^{n+1}}))^{\text{alg.}}$.
It can happen, under the hypotheses of Theorem \ref{mahler_nishioka}, that $f(\alpha)$ is well defined and algebraic for certain $\alpha\in\QQ^{\text{alg.}}\setminus\{0\}$.
For example, the formal series $$f(x)=\prod_{i=0}^\infty(1-2x^{2^i})\in\ZZ[[x]],$$ converging for $x\in\CC$ such that $|x|<1$
and satisfying the functional equation \begin{equation}f(x^2)=\frac{f(x)}{1-2x},\label{functionaleqf2}\end{equation} vanishes at every $\alpha$ such that
$\alpha^{2^i}=1/2$, $i\geq 0$. In particular, $f$ being non-constant and having infinitely many zeroes, it is transcendental. By Theorem \ref{mahler_nishioka}, $f(1/4)=\lim_{x\rightarrow 1/2}\frac{f(x)}{1-2x}$ is transcendental.

\begin{Remarque}\label{remarquenishioka}{\em 
Nishioka strengthened Theorem \ref{mahler_nishioka} allowing the rational function $R$ satisfying only the relaxed condition $h_Y(R)<d^2$ (see \cite[Theorem 1.5.1]{Ni} for an even stronger result).
The proof, more involved than the proof of Theorem \ref{mahler_nishioka}, follows most of the principles of it, with the following notable difference.
To achieve the proof, a more careful choice of the polynomials $P_N$ is needed.
In step (AP) it is again needed  to choose a sequence of polynomials $(P_N)_N$ with $P_N\in \QQ[X,Y]$ of degree $\leq N$
in $X$ and $Y$, such that the function $F_N(x)=P_N(x,f(x))$ vanishes at $x=0$ with order of vanishing $\geq c_1N^2$ for a constant $c_1$
depending on $\alpha$ and $f$. Since for $h_Y(R)\geq d$
the size of the coefficients of $P_N$ influences the conclusion, the use of Siegel's Lemma is now needed to accomplish this choice
\cite[Section 1.3]{Wa2}. To make good use of these refinements we need an improvement of the step (NV), since an explicit upper bound 
like $c_2N\log N$ for the integer $k$
such that $F(x^{d^s})=0$ for $s=0,\ldots,k$ is required (\footnote{This is not difficult to obtain; see \ref{hauptmodul} below for a similar, but more difficult estimate.}). }
\end{Remarque}

\subsection{Some further discussions.}

In this subsection we discuss about some variants of Mahler's method and applications to modular functions (in \ref{hauptmodul}).
We end with \ref{digression}, where we quote a
criterion of transcendence by Corvaia and Zannier quite different from Mahler's method, since it can be obtained as a corollary of Schmidt's subspace theorem.

\subsubsection{The ``stephanese" Theorem\label{hauptmodul}}

We refer to \cite{Silv} for a precise description of the tools concerning elliptic curves and modular forms and functions, involved in this subsection.

Let 
$$J(q)=\frac{1}{q}+744+\sum_{i\geq 1}c_iq^i\in(1/q)\ZZ[[q]]$$ be the {\em $q$-expansion} of the classical {\em hauptmodul}
for $\mathbf{SL}_2(\ZZ)$, converging for $q\in\CC$ such that $0<|q|<1$.
The following theorem was proved in 1996; see \cite{Ba}:
\begin{Theoreme}[BarrŽ-Sirieix, Diaz, Gramain and Philibert]\label{stephanese}
For $q$ complex such that $0<|q|<1$, one at least of the two complex numbers $q,J(q)$ is transcendental.
\end{Theoreme}
This {\em stephanese} theorem (\footnote{Sometimes, this result is called {\em stephanese theorem} from the name of the city of Saint-Etienne, where the authors of this result currently live.}) furnished a positive answer to Mahler's conjecture on values of the modular $j$-invariant (see \cite{Ma4}).
Although we will not say much more about, we mention that a similar conjecture was independently formulated by Manin, for $p$-adic values of $J$ at algebraic $\alpha$'s,
as the series $J$ also converges in all punctured $p$-adic unit disks, for every prime $p$.
Manin's conjecture is proved in \cite{Ba} as well. Manin's conjecture is relevant for its connections with the values of $p$-adic $L$-functions
and its consequences on $p$-adic variants of Birch and Swinnerton Dyer conjecture.

Mahler's conjecture was motivated by the fact that the function $J$ satis\-fies the autonomous non-linear {\em modular equation} $\Phi_2(J(q),J(q^2))=0$,
where \begin{eqnarray*}\Phi_2(X,Y)&=&X^3+Y^3-X^2Y^2+1488XY(X+Y)-162000(X^2+Y^2)+\\ & &40773375XY+8748000000(X+Y)-157464000000000.\end{eqnarray*} Mahler hoped to apply some suitable generalisation of Theorem \ref{mahler_nishioka}. It is still unclear, at the time being, if this intuition is correct; we remark that Theorem
\ref{mahler_nishioka} does not apply here.

The proof of Theorem \ref{stephanese} relies on a {\em variant} of Mahler's method that we discuss now. We first recall from \cite{Silv} that
there exists a collection of {\em modular equations} 
$$\Phi_{n}(J(q),J(q^n))=0,\quad n>0,$$ with explicitly calculable polynomials $\Phi_n\in\ZZ[X,Y]$ for all $n$. The stephanese team 
make use of the full collection of polynomials $(\Phi_n)_{n>0}$ so let us briefly explain how these functional equations occur.

For $q$ complex such that $0<|q|<1$, $J(q)$ is the {\em modular invariant} of an elliptic curve analytically isomorphic to 
the complex torus $\CC^\times/q^{\ZZ}$; if $z\in\CC$ is such that $\Im(z)>0$ and $e^{2\pi\mathrm{i}z}=q$, then there also is a torus analytic 
isomorphism $\CC^\times/q^{\ZZ}\equiv\CC/(\ZZ+z\ZZ)$. Since the lattice $\ZZ+nz\ZZ$ can be embedded in the lattice $\ZZ+z\ZZ$,
the natural map $\CC^\times/q^{\ZZ}\rightarrow\CC^\times/q^{n\ZZ}$ amounts to a cyclic {\em isogeny} of the corresponding elliptic curves
which, being projective smooth curves, can be endowed with {\em Weierstrass models} $y^2=4x^3-g_2x-g_3$ connected by algebraic relations
independent on the choice of $z$. At the level of the modular invariants, these algebraic relations for $n$ varying are precisely the 
modular equations, necessarily autonomous, defined over $\ZZ$ as a simple Galois argument shows.

Assuming that for a given $q$ with $0<|q|<1$, $J(q)$ is algebraic, 
means that there exists an {\em elliptic curve} $E$ analytically isomorphic to the torus $\CC^\times/q^\ZZ$, which is definable 
over a number field (it has Weierstrass model $y^2=4x^3-g_2x-g_3$ with $g_2,g_3\in\QQ^{\text{alg.}}$).
The discussion above, with the fact that the modular polynomials $\Phi_n$ are defined over $\ZZ$, implies that 
$J(q^n)$ is algebraic as well.

Arithmetic estimates involved in the (LB) step of the proof of Theorem \ref{stephanese} require a precise control, for $J(q)$ algebraic,
of the height of $J(q^n)$ and the degree $d_n$ of $\QQ(J(q^n),J(q))$ over $\QQ(J(q))$. the degree $d_n$ can be easily computed
counting lines in $\FF_p^2$ for $p$ prime dividing $n$; it thus divides the number $\psi(n)=\prod_{p|n}(1+1/p)$ and is bounded
from above by $c_3n^{1+\epsilon}$, for all $\epsilon>0$. As for the height $h_n=h(J(q^n))$, we said that the modular polynomial
$\Phi_n$ is related to a family of cyclic isogenies of degree $n$ connecting two families of elliptic curves. We then have, associated to the 
algebraic modular invariants $J(q),J(q^n)$, two isogenous elliptic curves defined over a number field, and the isogeny has degree
$n$. Faltings theorem asserting that the modular heights of two isogenous elliptic curves may 
differ of at most the half of the logarithm of a minimal degree of isogeny gives the bound $c_3(h(J(q))+(1/2)\log n)$ for the logarithmic height 
$h(J(q^n))$ (this implies the delicate estimates the authors do in \cite{Ba}).

With these information in mind, the proof of Theorem \ref{stephanese} proceeds as follows. As in remark \ref{remarquenishioka},
we use standard estimates of the growth of the absolute values of the (integral) coefficients of the $q$-expansions
of the normalised Eisenstein series $E_4,E_6$ of weights $4,6$ to
apply Siegel's Lemma and construct a sequence of auxiliary polynomials (AP) $(P_N)_{N\geq 1}$ in $\ZZ[X,Y]\setminus\{0\}$
with $\deg_XP_N,\deg_YP_N\leq N$, such that $F_N(x):=P_N(x,xJ(x))$ vanishes with order $\geq N^2/2$ at $x=0$.

The (UB) estimate is then exactly as in the Proof of Theorem \ref{mahler_nishioka}. All the authors of \cite{Ba}
need to achieve their proof is the (NV) step; and it is here that a new idea occurs.
They use that the coefficients of $J$ are {\em rational integers} to 
deduce a sharp estimate of the biggest integer $n$ such that $F_N(x)$ 
vanishes at $q^m$ for all $m=1,\ldots,n-1$. This idea, very simple and appealing to Schwarz lemma,
does not seem to occur elsewhere in Mahler's theory; it was later generalised by Nesterenko in the proof of his famous theorem in \cite{Nest, Nesterenko:Introduction3}, implying
the algebraic independence of the three numbers $\pi,e^{\pi},\Gamma(1/4)$ and the stephanese theorem
(\footnote{We take the opportunity to notice that a proof of an analog of the stephanese theorem for the so-called ``Drinfeld modular invariant" 
by Ably, Recher and Denis is contained in \cite{Ably}.}). We will come back to the latter result in Section \ref{algindep}.

\subsubsection{Effects of Schmidt's Subspace Theorem\label{digression}}

We mention the following result in \cite{CoZa} whose authors Corvaja and Zannier deduce from 
Schmidt's Subspace Theorem.

\begin{Theoreme}[Corvaja and Zannier]\label{corvajazannier}
Let us consider a formal series $f\in\QQ^{\text{alg.}}((x))\setminus\QQ^{\text{alg.}}[x,x^{-1}]$ 
and assume that $f$ converges for $x$ such that $0<|x|<1$.
Let $L\subset \CC$ be a number field and $S$ a finite set of places of $L$ containing the archimedean ones.
Let $\mathcal{A}\subset\NN$ be an infinite subset.
Assume that:
\begin{enumerate}
\item $\alpha\in L$, $0<|\alpha|<1$
\item $f(\alpha^n)\in L$ is an $S$-integer for all $n\in\mathcal{A}$.
\end{enumerate}
Then, 
$$\liminf_{n\in\mathcal{A}}\frac{h(f(\alpha^n))}{n}=\infty$$
\end{Theoreme}

This theorem has as an immediate application with $\mathcal{A}=\{d,d^2,d^3,\ldots\}$, $d>1$ being an integer.
If $f\in\QQ^{\text{alg.}}[[x]]$ is not a polynomial, converges for $|x|<1$ and is such that
$f(x^d)=R(x,f(x))$ with $R\in\QQ^{\text{alg.}}(X,Y)$
with $h_Y(R)<d$, then, $f(\alpha^{d^n})$ is transcendental for $\alpha$ algebraic with $0<|\alpha|<1$ and for all
$n$ big enough. This implies a result (at least apparently) stronger than Theorem \ref{mahler_nishioka}; indeed, the hypothesis
that the coefficients of the series $f$ all lie in a given number field is dropped. 



\section{Transcendence theory in positive characteristic\label{anotherexample}}

The reduction modulo $2$ in $\FF_2[[x]]$ of the formal series $f_{\text{TM}}(x)\in\ZZ[[x]]$ is an algebraic formal series.
In this section we will see that several interesting transcendental series in positive characteristic are analogues 
of the series satisfying the functional equation (\ref{functionaleqf2}).

Let $q=p^e$ be an integer power of a
prime number $p$ with $e>0$, let $\FF_q$ be the field with $q$ elements.
Let us write $A=\FF_q[\theta]$ and $K=\FF_q(\theta)$, with $\theta$
an indeterminate over $\FF_q$, and define an absolute value
$|\cdot|$ on $K$ by $|a|=q^{\deg_\theta a}$, $a$ being in $K$, so that
$|\theta| = q$.  Let $K_\infty :=\FF_q((1/\theta))$ be the
completion of $K$ for this absolute value, let $K_\infty^{\text{\tiny alg.}}$ be an
algebraic closure of $K_\infty$, let $C$ be the completion of
$K_\infty^{\text{\tiny alg.}}$ for the unique extension of $|\cdot|$ to $K_\infty^{\text{\tiny alg.}}$, and let $K^{\text{\tiny alg.}}$ be the algebraic closure of $K$
embedded in $C$.
There is a unique degree map $\deg_\theta:C^\times\rightarrow\QQ$ which extends 
the map $\deg_\theta:K^\times\rightarrow\ZZ$.

Let us consider the power series
$$f_{\text{De}}(x)=\prod_{n=1}^\infty(1-\theta x^{q^n}),$$
which converges for all $x\in C$ such that $|x|<1$ and satisfies, just as in (\ref{functionaleqf2}), the functional equation:
\begin{equation}f_{\text{De}}(x^q)=\frac{f_{\text{De}}(x)}{1-\theta x^q}\label{functionalequationl}\end{equation}
(the subscript De stands for Denis, who first used this series for transcendence purposes).

For $q=2$, we notice that $f_{\text{De}}(x)=\sum_{n=0}^\infty \theta^{b_n}x^{2n},$ where $$(b_n)_{n\geq 0}=0, 1, 1, 2, 1, 2, 2, 3, 1, 2, 2, 3,\ldots$$
is the sequence with $b_n$ equal to the sum of the digits of the binary expansion of $n$ (and whose reduction modulo $2$ precisely is Thue-Morse sequence of Section \ref{firstexample}). It is very easy to show that $f_{\text{De}}$ is transcendental, because it is plain that $f_{\text{De}}$ has infinitely many zeros $\theta^{-1/q},\theta^{-1/q^2},\ldots$ (just as the function occurring at the end of \ref{moregeneralres}). We shall prove:

\begin{Theoreme} For all $\alpha\in K^{\text{alg.}}$ with $0<|\alpha|<1$, $f_{\text{De}}(\alpha)$ is transcendental.\label{secondexample}\end{Theoreme}
 
 \subsection{Proof of Theorem \ref{secondexample}}
 
 The proof of Theorem \ref{secondexample} follows the essential lines of Section \ref{firstexample}, once the necessary tools are introduced.
 
 \subsubsection{Transcendence of functions\label{transcendencefunctionsvartheta}}
 
Not all the arguments of \ref{ftranscendental} work well to show the transcendence of formal series such as
$f_{\text{De}}$; in particular, the so-called Riemann-Hurwitz-Hasse formula does not give much information for functional equations
 such as $f_{\text{De}}(x^d)=af_{\text{De}}(x)+b$ with the characteristic that divides $d$.
Since in general it is hard to detect zeros of Mahler's functions,
 we report another way to check the transcendency of $f_{\text{De}}$, somewhat making use of ``automatic methods",
 which can also be generalised as it does not depend on the location of the zeroes. To simplify the presentation,
 we assume, in the following discussion, that $q=2$ but at the same time, we relax certain conditions so that, 
 in all this subsection, we denote by $\vartheta$ an element of $C$ and 
 by $f_\vartheta$ the formal series 
 $$f_\vartheta(x)=\prod_{n=0}^\infty(1-\vartheta x^{q^n})=\sum_{n=0}^\infty\vartheta^{b_n}x^n\in F[[x]]\subset C[[x]]$$ with 
 $F$ the perfect field $\cup_{n\geq 0}\FF_2(\vartheta^{1/2^n})$, converging for $x\in C$ with $|x|<1$.
 In particular, we have
  $$f_\theta(x)=f_{\text{De}}(x^{1/2})\in\FF_2[\theta][[x]].$$
  We shall prove:
  \begin{Theoreme}
  The formal series $f_\vartheta$ is algebraic over $F(x)$ if and only if $\vartheta$ belongs to $\overline{\FF_q}$,
  embedded in $C$. 
  \label{alasharif}\end{Theoreme}
 \noindent\emph{Proof.} If $\vartheta\in\overline{\FF_q}$, it is easy to show that $f_\vartheta$ is algebraic,
 so let us assume by contradiction that $f_\vartheta$ is algebraic, with $\vartheta$ that belongs to $C\setminus\overline{\FF_q}$.
 
 We have the functional equation:
 \begin{equation}(1-\vartheta x)f_\vartheta(x^2)=f_\vartheta(x).\label{newfunctional}\end{equation}
 
 We introduce the operators $$f=\sum_{i}c_ix^i\in F((x))\mapsto f^{(k)}=\sum_{i}c_i^{2^k}x^i\in F((x)),$$
 well defined for all $k\in\ZZ$. Since $f(x^2)=f^{(-1)}(x)^2$ for any $f\in F[[x]]$, we deduce from (\ref{newfunctional})
 the collection of functional equations
 \begin{equation}\label{v}
 f_\vartheta^{(-1-k)}(x)^2(1-\vartheta^{1/2^k}x)=f_\vartheta^{(-k)}(x),\quad k\geq0.\end{equation}
 
 For any $f\in F[[x]]$ there exist two series $f_0,f_1\in F[[x]]$, uniquely determined, with the property
 that
 $$f=f_0^2+xf_1^2.$$ We define $E_i(f):=f_i$ ($i=0,1$).
 It is plain that, for $f,g\in F[[x]]$,
 \begin{eqnarray*}
E_i(f+g)&=&E_i(f)+E_i(g),\quad (i=0,1),\\
E_0(fg)&=&E_0(f)E_0(g)+xE_1(f)E_1(g),\\
E_1(fg)&=&E_0(f)E_1(g)+E_1(f)E_0(g),\\
E_0(f^2)&=&f,\\
E_1(f^2)&=&0.
\end{eqnarray*}
Therefore, 
$$E_i(f^2g)=fE_i(g),\quad i=0,1.$$
 
By (\ref{v}) we get
$$\begin{array}{l}E_0(f_\vartheta^{(-k)})=f_\vartheta^{(-1-k)}E_0(1-\vartheta^{1/2^k}x)=f_\vartheta^{(-1-k)},\\ E_1(f_\vartheta^{(-k)})=E_1(1-\vartheta^{1/2^k}x)f_\vartheta^{(-1-k)}=-\vartheta^{1/2^k}f_\vartheta^{(-1-k)},\end{array}$$
and we see that if $V$ is a $F$-subvector space of $F[[x]]$ containing $f_\vartheta$ and 
stable under the action of the operators $E_0,E_1$, then $V$ contains the $F$-subvector space
generated by $f_\vartheta,f_\vartheta^{(-1)},f_\vartheta^{(-2)},\ldots$.

By a criterion
for algebraicity of Sharif and Woodcock \cite[Theorem 5.3]{ShaWood}
there is a subvector space $V$ as above, with finite dimension, containing $f_\vartheta$. The formal series $f_\vartheta^{(-k)}$ are
$F$-linearly dependent and there exists $s>0$ such that $f_\vartheta,f_\vartheta^{(1)},\ldots,f_\vartheta^{(s-1)}$ are $F$-linearly
dependent.

Going back to the explicit $x$-expansion of $f_\vartheta$, the latter condition is equivalent to 
the existence of $c_0,\ldots,c_{s-1}\in F$, not all zero, such that for all $n\geq 0$:
$$\sum_{i=0}^{s-1}c_i\vartheta^{2^ib_n}=0.$$

The sequence $b:\NN\cup\{0\}\rightarrow\NN\cup\{0\}$ is known to be surjective, so that
the Moore determinant
$$\det((\vartheta^{2^ij}))_{0\leq i,j\leq s-1}$$
vanishes. But this means that $1,\vartheta,\vartheta^2,\ldots,\vartheta^{s-1}$ are $\FF_2$-linearly 
dependent (Goss, \cite[Lemma 1.3.3]{Go}), or in other words, that $\vartheta$ is algebraic over $\FF_2$; a contradiction which completes
the proof that $f_\vartheta$ and in particular $f_{\text{De}}$ are transcendental over $F(x)$ (and the fact that the image of $b$ has
infinitely many elements suffices to achieve the proof).

\subsubsection{Heights under a more general point of view\label{logarithmicheight}}

A good framework to generalise logarithmic heights to other base fields is that described by Lang in \cite[Chapter 3]{La} and by Artin and Whaples \cite[Axioms 1, 2]{AW1}. Let $\mathcal{K}$ be any field together with a
proper set of non-equivalent places $M_{\mathcal{K}}$. Let us choose, for every place $v\in M_{\mathcal{K}}$ an absolute value $|\cdot|_v\in v$
and assume that for all $x\in\mathcal{K}^\times$, the following product formula holds (cf. \cite{La} p. 23):
\begin{equation}\prod_{v\in M_{\mathcal{K}}}|x|_v=1,\quad x\in \mathcal{K}^\times,\label{eq:productformula}
\end{equation}
with the additional property that if $\alpha$ is in $\mathcal{K}^\times$, then 
$|\alpha|_v=1$ for all but finitely many $v\in M_\mathcal{K}$. Let us suppose that $M_\mathcal{K}$ contains
at least one absolute value associated to either a discrete, or an archimedean valuation of $\mathcal{K}$.
It is well known that under these circumstances \cite{AW1}, $\mathcal{K}$ is either 
a number field, or a function field of one variable over a field of constants.

Given a finite extension $\mathcal{L}$
of $\mathcal{K}$, there is a proper set $M_{\mathcal{L}}$ of absolute values on $\mathcal{L}$, extending those of $M_{\mathcal{K}}$, again satisfying the product formula
\begin{equation}\label{ultrametricproductformula}
\prod_{v\in M_L}|\alpha|_v^{d_v}=1,
\end{equation}
where,
 if $v$ is the place of $\mathcal{K}$ such that $w|_\mathcal{K}=v$ (one then writes $w|v$), we have defined $d_w=[L_w,K_v]$, so that $\sum_{w|v}d_v=[\mathcal{L}:\mathcal{K}]$.
 
An analogue of the absolute logarithmic height $h$ is available, by the following definition (see \cite[Chapter 3]{La}). Let $(\alpha_0:\cdots:\alpha_n)$ be a projective point defined over
$\mathcal{L}$. Then we define:
$$h(\alpha_0:\cdots:\alpha_n)=\frac{1}{[\mathcal{L}:\mathcal{K}]}\sum_{w\in M_{\mathcal{L}}}d_w\log\max\{|\alpha_0|_w,\ldots,|\alpha_n|_w\}.$$

Again, we have a certain collection of properties making this function useful in almost every proof of transcendence over function fields.

First of all, product formula 
(\ref{ultrametricproductformula})
implies that $h(\alpha_0:\cdots:\alpha_n)$ does not depend on the choice of the field $\mathcal{L}$ and defines
a map $$h:\PP_n(\mathcal{K}^{\text{\tiny alg.}})\rightarrow\RR_{\geq 0}.$$

We write $h(\alpha):=h(1:\alpha)$.
If the absolute values of $M_\mathcal{K}$ are all ultrametric, it is easy to prove, 
with the same indications as in \ref{height}, that
for $\alpha,\beta\in \mathcal{K}^{\text{alg.}}{}^\times$:
\begin{eqnarray*}
h(\alpha+\beta),h(\alpha\beta)&\leq &h(\alpha)+h(\beta),\\
h(\alpha^n)&=&|n|h(\alpha),\quad n\in\ZZ.
\end{eqnarray*}
More generally, if $P$ is a polynomial in $\mathcal{L}[X_1,\ldots,X_n]$ and if $(\alpha_1,\ldots,\alpha_n)$ is a point of $\mathcal{L}^n$, we write $h(P)$ for the height of the projective 
point whose coordinates are $1$ and its coefficients. We have:
\begin{equation}\label{heightglobal}
h(P(\alpha_1,\ldots,\alpha_n))\leq h(P)+\sum_{i=1}^n(\deg_{X_i}P)h(\alpha_i).
\end{equation}

Product formula (\ref{ultrametricproductformula}) also provides a Liouville's type inequality.
Let $[\mathcal{L}:\mathcal{K}]_{\text{sep.}}$ be the separable degree of $\mathcal{L}$ over $\mathcal{K}$.
 Let us choose a distinguished absolute value 
$|\cdot|$ of $\mathcal{L}$ and $\beta\in L^\times$. We have:
\begin{equation}\label{liouville2}
\log|\beta|\geq -[\mathcal{L}:\mathcal{K}]_{\text{sep.}}h(\beta).
\end{equation}
The reason of the presence of the separable degree in (\ref{liouville2}) is the following. If $\alpha\in\mathcal{L}^\times$ is separable
over $\mathcal{K}$ then $\log|\alpha|\geq -[\mathcal{K}(\alpha):\mathcal{K}]h(\alpha)= -[\mathcal{K}(\alpha):\mathcal{K}]_{\text{sep.}}h(\alpha)$.
Let $\beta$ be any element of $\mathcal{L}^\times$. There exists $s\geq 0$ minimal with $\alpha=\beta^{p^s}$ separable and we get
$p^s\log|\beta|\geq-[\mathcal{K}(\alpha):\mathcal{K}]h(\alpha)=-p^s[\mathcal{K}(\beta):\mathcal{K}]_{\text{sep.}}h(\beta)$.

\subsubsection{Transcendence of the values of $f_{\text{De}}$\label{transcendencel}}

We now follow Denis and
we take $\mathcal{K}=K$, $M_K$ the set of all the places of $K$ and we choose in each of these places an absolute value normalised so that product formula (\ref{eq:productformula}) holds, with the distinguished absolute value $|\cdot|$ chosen so that $|\alpha|=q^{\deg_\theta\alpha}$ for $\alpha\in K^\times$.

As we already did in \ref{proofofexample}, we choose
 for all $N\geq 0$, a polynomial $P_N\in K[X,Y]$, non-zero, of degree $\leq N$
in both $X,Y$, such that the order of vanishing $\nu(N)<\infty$ of the function $F_N(x):=P_N(x,f_{\text{De}}(x))$
at $x=0$ satisfies $\nu(N)\geq N^2$. We know that this is possible by simple linear algebra arguments as 
we did before.

Let $\alpha\in K^{\text{\tiny alg.}}$ be such that $0<|\alpha|<1$; as
in \ref{proofofexample}, the sequence $(P_N)_N$ need not to depend on it but the choice of $N$ we will do does.

By the identity principle of analytic functions on $C$, if $\epsilon$ is a positive real number, for $n$ big enough depending on $N$ and $\alpha,l,\epsilon$, we have
$F_N(\alpha^{q^{n+1}})\neq0$ and $$\log|F_N(\alpha^{q^{n+1}})|\leq \nu(N)q^{n+1}\log|\alpha|+\log|c_{\nu(N)}|+\epsilon,$$ where $c_{\nu(N)}$ is a non-zero element of $K$ depending on $N$ (it is the leading coefficient of the formal series $F_N$).

Let us assume by contradiction that 
$$f_{\text{De}}(\alpha)\in K^{\text{\tiny alg.}},$$ let $L$ be a finite extension of $K$ containing $\alpha$ and $f_{\text{De}}(\alpha)$.

By the variant of Liouville's inequality (\ref{liouville2}) and from the basic facts on the height $h$ explained above
\begin{eqnarray*}
\lefteqn{\log|F_N(\alpha^{q^{n+1}})|\geq}\\
&\geq &-[L:K]_{\text{sep.}}\left(\deg_\theta P_N+Nh(\alpha^{q^{n+1}})+Nh\left(\frac{f_{\text{De}}(\alpha)}{(1-\theta\alpha^q)\cdots(1-\theta\alpha^{q^{n+1}})}\right)\right).
\end{eqnarray*}
Dividing by $Nq^{n+1}$ and using that $\nu(N)\geq N^2$ we get, for all $n$ big enough,
\begin{eqnarray*}
\lefteqn{N\log|\alpha|+\log|c_{\nu(N)}|+\epsilon\geq}\\ & \geq & -[L:K]_{\text{sep.}}(\deg_\theta P_NN^{-1}q^{-n-1}+(1+(q-q^{-n-1})/(q-1))h(\alpha)+\\ & & q^{-n-1}h(f_{\text{De}}(\alpha))+(n+1)N^{-1}q^{-n-1}h(\theta)).\end{eqnarray*}
 Letting $n$ tend to infinity, we obtain the inequality:
 $$N\log|\alpha|\geq -[L:K]_{\text{sep.}}\left(1+\frac{q}{q-1}\right)h(\alpha)$$ for all $N> 0$. Just as in the proof of Theorem \ref{firstexamplebis}, if $N$ is big enough, this is contradictory with 
 the assumptions showing that $f_{\text{De}}(\alpha)$ is transcendental.

\subsubsection{A first application to periods}

The transcendence of values of $f_{\text{De}}$ at algebraic series has interesting applications, especially when one looks at what happens with the base point $\alpha=\theta^{-1}$. Indeed, let
 \begin{equation}\widetilde{\pi}= \theta(-\theta)^{1/(q-1)} \prod_{i=1}^\infty
 ( 1 - \theta^{1-q^i})^{-1}\label{prodottopi}\end{equation} be a fundamental period of Carlitz's module (it is defined up to multiplication by an element of $\FF_q^\times$).
Then, $$
\widetilde{\pi}=\theta(-\theta)^{1/(q-1)}f_{\text{De}}(\theta^{-1})^{-1},$$
  so that it is transcendental over $K$.
 
If $\alpha=\theta^{-1}$, $h(\alpha)=\log q$ so that to show that $\widetilde{\pi}$ is not in $K$, it suffices to choose $N\geq 4$ if 
$q=2$ and $N\geq 3$ if $q\neq2$ in the proof above. Let us look, for $q\neq2$ given, at a polynomial
 (depending on $q$) $P\in A[X,Y]\setminus\{0\}$ with relatively prime coefficients 
in $X$ of degree $\leq 3$ in $X$ and in $Y$, such that $P(u,f_{\text{De}}(u))$ vanishes at $u=0$ with the biggest possible order $\nu> 9$ (which also depends on $q$).

It is possible to prove that for all $q\geq 4$,
$$P=X^3(Y-1)^3\in\FF_q[X,Y].$$ 
This means that to show that $f_{\text{De}}(\theta^{-1})\not\in K$,
it suffices to work with the polynomial $Q=Y-1$ (\footnote{Other reasons, related to the theory of Carlitz module, allow to show directly that $\widetilde{\pi}\not\in \FF_q((\theta^{-1}))$ for $q>2$}). Indeed, $f_{\text{De}}(u)-1=-\theta u^q+\cdots$.
Therefore, for $n$ big enough, if by contradiction $f_{\text{De}}(\theta^{-1})\in K$, then $\log|f_{\text{De}}(\alpha^{q^{n+1}})-1|\geq -(1+q/(q-1))h(\alpha^{q^{n+1}})$ which is contradictory even taking $q=3$, but not for $q=2$, case that we skip. 

If $q=3$, we get a completely different kind of polynomial $P$ of degree $\leq 3$ in each indeterminate:
$$P=2+\theta^2+\theta X^3+2\theta^3X^3+2\theta^2Y+Y^3+2\theta X^3Y^3.$$
It turns out that
$$P(u,f_{\text{De}}(u))=(\theta^6-\theta^4)u^{\nu}+\cdots$$
with $\nu=36$ so that the order of vanishing is three times as big as the quantity expected from the computations with $q\geq 4$: $12=3\cdot 3+3$.
How big can $\nu$ be? It turns out that this question is important, notably in the search for quantitative measures of
transcendence and algebraic independence; we will discuss about this problem in Section \ref{criteria}.

\subsection{A second transcendence result\label{foursteps}}

With essentially the same arguments of Section \ref{anotherexample}, it is possible to 
deal with a more general situation and prove the Theorem below. We first explain the data we will work with.

Let $F$ be a field and $t$ an indeterminate.
We denote by $F\langle\langle t\rangle\rangle$ the field of {\em Hahn generalised series}. This is the set of formal series
$$\sum_{i\in S}c_it^i,\quad c_i\in F,$$ with $S$ a {\em well ordered subset} of $\QQ$
(\footnote{By definition, every nonempty subset of $S$
has a least element for the order $\leq$.}), endowed with the standard addition and Cauchy's multiplication
from which it is plain that every non-zero formal series is invertible.

We have a field $\FF_q\langle\langle t\rangle\rangle$-automorphism
$$\tau:C\langle\langle t\rangle\rangle\rightarrow C\langle\langle t\rangle\rangle$$
defined by
$$\alpha=\sum_{i\in S}c_it^i\mapsto\tau\alpha=\sum_{i\in S}c_i^qt^i.$$

Assume that, with the notations of \ref{logarithmicheight}, $\mathcal{K}=C(t)$, with $t$ an independent
indeterminate.
%
Let $|\cdot|$ be an absolute value associated to the $t$-adic valuation
and $\widehat{\mathcal{K}}$ the completion of $\mathcal{K}$ for this absolute value.
Let $\KK$ be the completion of an algebraic closure of $\widehat{\mathcal{K}}$ for the extension 
of $|\cdot|$, so that we have an embedding of $\mathcal{K}^{\text{alg.}}$ in $\KK$. 
We have an embedding $\iota:\KK\rightarrow C\langle\langle t\rangle\rangle$ (see Kedlaya, \cite[Theorem 1]{Kedlaya}); there
exists a rational number $c>1$ such that 
if $\alpha$ is in $\KK$ and $\iota(\alpha)=\sum_{i\in S}c_it^i$, then $|\alpha|=c^{-i_0}$, where
$i_0=\min (S)$. 
We identify $\KK$
with its image by $\iota$. It can be proved that $\tau\KK\subset\KK$, $\tau\mathcal{K}^{\text{alg.}}\subset\mathcal{K}^{\text{alg.}}$
and $\tau\mathcal{K}\subset\mathcal{K}$.

The definition of $\tau$ implies immediately that, for all $\alpha\in\KK$,
\begin{equation}\label{propertytau1}|\tau\alpha|=|\alpha|.\end{equation}

We choose $M_\mathcal{K}$ a complete set of non-equivalent 
absolute values of $\mathcal{K}$ such that the product formula (\ref{eq:productformula}) holds.
On $\PP_n(\mathcal{K}^{\text{alg.}})$, we have the absolute logarithmic height whose 
main properties have been described in \ref{logarithmicheight}.

There is a useful expression for the height $h(\alpha)$ of a non-zero element $\alpha$
in $\mathcal{K}^{\text{alg.}}$ of degree $D$. If $P=a_0X^d+a_1X^{d-1}+\cdots+a_{d-1}X+a_d$ is 
a polynomial in $C[t][X]$ with relatively prime coefficients such that $P(\alpha)=0$, we have:
\begin{equation}\label{formulaheights}
h(\alpha)=\frac{1}{D}\left(\log|a_0|+\sum_{\sigma:\mathcal{K}^{\text{alg.}}\rightarrow\KK}\log\max\{1,|\sigma(\alpha)|\}\right),
\end{equation}
where the sum runs over all the $\mathcal{K}$-embeddings of $\mathcal{K}^{\text{alg.}}$ in $\KK$.
The proof of this formula follows the same ideas as that of \cite[Lemma 3.10]{Wa}.

Let $\alpha$ be in $\mathcal{K}^\text{alg.}$.
From (\ref{propertytau1}) and (\ref{formulaheights}), it follows that:
\begin{equation}\label{propertytau}
h(\tau\alpha)=h(\alpha).
\end{equation}

Let us also consider, over the ring of formal series
$\KK[[x]]$, the $\FF_q\langle\langle t\rangle\rangle$-extension of $\tau$ defined in the following way:
$$f:=\sum_ic_ix^i\mapsto \tau f:=\sum_i(\tau c_i)x^{qi}.$$
We can now state the main result of this section.

\begin{Theoreme}\label{theo1} Let $f\in \mathcal{K}[[x]]$ be converging for $x\in \KK$, $|x|< 1$,
let $\alpha\in \mathcal{K}$ be such that $0<|\alpha|<1$.
Assume that:
 \begin{enumerate}
 \item $f$ is transcendental over $\mathcal{K}(x)$,
 \item $\tau f=af+b$, where $a,b$ are elements of $\mathcal{K}(x)$.
 \end{enumerate}
 Then, for all $n$ big enough,
 $(\tau^nf)(\alpha)$ is transcendental over $\mathcal{K}$.
\end{Theoreme}
\noindent\emph{Proof.} We begin with a preliminary discussion about heights.
Let $r=r_0x^n+\cdots+r_n$ be a polynomial in $\mathcal{K}[x]$. We have, for all 
$j\geq 0$, $\tau^jr=(\tau^jr_0)x^{q^jn}+\cdots+(\tau^jr_n)$. Therefore, if $\alpha$ is an
element of $\mathcal{K}^{\text{alg.}}$, we deduce from (\ref{heightglobal}), (\ref{propertytau}) and from 
elementary height estimates:
\begin{eqnarray*}
h((\tau^jr)(\alpha))&\leq &h(1:\tau^jr_0:\cdots:\tau^jr_n)+q^jnh(\alpha),\\
&\leq&\sum_{i=0}^nh(\tau^jr_i)+q^jnh(\alpha)\\
&\leq&\sum_{i=0}^nh(r_i)+q^jnh(\alpha),
\end{eqnarray*}
where we have applied 
(\ref{propertytau}). Therefore, if $a$ is a rational function in $\mathcal{K}(x)$ such that $(\tau^ja)(\alpha)$
is well defined, we have
\begin{equation}\label{intermediateheightestimate}
h((\tau^ja)(\alpha))\leq c_1+q^jc_2,
\end{equation}
where $c_1,c_2$ are two constants depending on $a,\alpha$ only.

The condition on $f$ implies that, for all $k\geq 0$,
\begin{equation}\label{formulatauk}
\tau^kf=f\prod_{i=0}^{k-1}(\tau^i a)+\sum_{i=0}^{k-1}(\tau^ib)\prod_{j=i+1}^{k-1}(\tau^ja)
\end{equation}
(where empty sums are equal to zero and empty products are equal to one). Hence,
the field $$\mathcal{L}=\mathcal{K}(\alpha,f(\alpha),(\tau f)(\alpha),(\tau^2 f)(\alpha),\ldots)$$
is equal to $\mathcal{K}(\alpha,(\tau^n f)(\alpha))$ for all $n$ big enough. 

The transcendence of $f$ implies that  $a\neq0$. If $\alpha$ is a zero or a pole of $\tau^ka$ and a pole of $\tau^kb$
for all $k$, then it is a simple exercise left to the reader to prove that $|\alpha|=1$, case that we have excluded.

Let us suppose by contradiction that the conclusion of the theorem is false. Then, $\alpha$ is not a pole or a zero of
$\tau^na,\tau^nb$, 
$(\tau^nf)(\alpha)$ is algebraic over $\mathcal{K}$ for all $n$ big enough, and $\mathcal{L}$ is a finite extension of $\mathcal{K}$.

An estimate for the height
of this series can be obtained as follows.

A joint application of (\ref{propertytau}), (\ref{intermediateheightestimate}) and (\ref{formulatauk}) yields:
\begin{eqnarray*}
h((\tau^kf)(\alpha))&\leq&h(f(\alpha))+\sum_{i=0}^{k-1}h((\tau^ia)(\alpha))+\sum_{i=0}^{k-1}h((\tau^ib)(\alpha))\\
&\leq&c_3+c_4k+c_5q^k.
\end{eqnarray*}
Therefore, if $P$ is a polynomial in $\mathcal{K}[X,Y]$ of degree $\leq N$ in $X$ and $Y$,
writing $F_k$ for the formal series $\tau^kP(x,f(x))=P^{\tau^k}(x^{q^k},(\tau^kf)(x))$ ($P^{\tau^k}$
is the polynomial obtained from $P$, replacing the coefficients by their images under $\tau^k$),
we get:
\begin{eqnarray}
h(F_k(\alpha))&\leq&h(P^{\tau^k})+q^k(\deg_XP)h(\alpha)+(\deg_YP)h((\tau^kf)(\alpha))\nonumber\\
&\leq&c_7(P)+c_6Nq^k,\label{estimatePtauk}
\end{eqnarray}
where $c_7$ is a constant depending on $P$.

Let $N$ be a positive integer. There exists $P_N\in\mathcal{K}[X,Y]$ with 
partial degrees in $X,Y$ not bigger than $N$, with the additional property that
$F_N(x):=P_N(x,f(x))=c_{\nu(N)}x^{\nu(N)}+\cdots$, with $\nu(N)\geq N^2$ and $c_{\nu(N)}\neq0$.

Let us write $F(x)=\sum_{i\geq 0}c_ix^i$. 
In ultrametric analysis, Newton polygons suffice to locate the absolute values of the zeroes of
Taylor series.
The Newton polygons of the series $\sum_{i\geq 0}(\tau^kc_i)x^i\in\mathcal{K}[[x]]$ for $k\geq 0$ are all equal by (\ref{propertytau1}). By \cite[Propositions 2.9, 2.11]{Go},
we have $\sum_{i\geq 0}(\tau^kc_i)\alpha^{q^ki}\neq0$ for $k$ big enough.
Now, since for $k\geq 0$, $$(\tau^kF_N)(x)=(\tau^kc_{\nu(N)})x^{\nu(N)q^k}+\cdots,$$ we find, 
when the logarithm is well defined and by (\ref{propertytau}),
that 
\begin{equation}\label{equ1}
-\infty<\log|(\tau^k F_N)(\alpha)|\leq \nu(N)q^{k}\log|\alpha|+\log|c_{\nu(N)}|+\epsilon.
\end{equation}
On the other side, by (\ref{estimatePtauk}), 
\begin{equation}\label{equ2}
h((\tau^kF_N)(\alpha))\leq c_8(N)+c_9Nq^k,
\end{equation}
where $c_8$ is a constant depending on $f,\alpha,N$ and $c_9$ is a constant depending on $f,\alpha$.

A good choice of $N$ (big) and inequality (\ref{liouville2}) with $k$ big enough depending on $N$ give a contradiction
(\footnote{It is likely that Nishioka's proof of Theorem \ref{mahler_nishioka} can be adapted to strengthen Theorem
\ref{theo1}, but we did not enter into the details
of this verification.}).

\subsubsection{Applications of Theorem \ref{theo1}.}

We look at solutions $f\in \mathcal{K}[[x]]$ of $\tau$-difference equations
\begin{equation}
\tau f=af+b,\quad a,b\in \mathcal{K}(x)\label{abp2}.\end{equation} Theorem \ref{theo1} allows to 
give some information about the arithmetic properties of their values.

\medskip

\noindent\emph{First application.} Assume that in (\ref{abp2}), $a,b\in\FF_q(t)$. Then, since
$\FF_q(t)$ is contained in the field of constants of $\tau$, solutions of this difference equation
are related to the variant of Mahler's method of Section \ref{anotherexample}.

If $a=(1-t^{-1}x)^{-1},b=0$, the equation above has the solution
$$f_{\text{De2}}(x)=\prod_{n=1}^\infty(1-t^{-1}x^{q^n}),$$ which converges for $x\in\KK$, $|x|<1$.

If $x=t$, Theorem \ref{theo1} yields the transcendence of $f_{\text{De2}}(t)=\prod_{n=1}^\infty(1-t^{q^n-1})\in\FF_q[[t]]$
over $\FF_q(t)$ and we get (again) the transcendence of $\widetilde{\pi}$ over $K$
(we also notice the result of the paper \cite{Allouche}, which allows some other applications).
More generally, all the examples of functions in \cite[Section 3.1]{Pellarin} have a connection with this example.

\medskip

\noindent\emph{Second application.} Theorem \ref{theo1} also has some application which does not seem to immediately
follow from results such as Theorem \ref{mahler_nishioka}. Consider equation (\ref{abp2}) with $b=0$ and $a=(1+\vartheta x)^{-1}$, where 
$\vartheta\in\FF_q(t,\theta)$ is non-zero. We have the following solution of (\ref{abp2}) in $\mathcal{K}[[x]]$:
$$\phi(x)=\prod_{n=0}^\infty(1+(\tau^n\vartheta) x^{q^n}).$$
It is easy to show that $\phi=\sum_{j\geq 0}c_jx^j$ is a formal series of $\mathcal{K}[[x]]$ converging for $x\in\KK$, $|x|<1$. 

The coefficients $c_j$ can be computed in the following way. We have $c_j=0$ if the $q$-ary expansion
of $j$ has its set of digits not contained in $\{0,1\}$.  Otherwise, if $j=j_0+j_1q+\cdots+j_nq^n$ with 
$j_0,\ldots,j_n\in\{0,1\}$, we have, writing $\vartheta_i$ for $\tau^ir$, $c_j=\vartheta_0^{j_0}\vartheta_1^{j_1}\cdots \vartheta_n^{j_n}$.
Therefore, if $\psi(x)=\sum_{k=0}^\infty \vartheta_0\vartheta_1\cdots \vartheta_k x^{1+q+\cdots+q^k}$,
we have $\phi(x)-1=\sum_{j=0}^\infty\tau^j\psi(x)$.

For $\vartheta=-t^{-1}(1+t/\theta)^{-1}$, the series $\phi(x)$ vanishes at every $x_n=t^{1/q^n}(1+t^{1/q^n}/\theta)$. The $x_n$'s are elements
of $\KK$ which are distinct with absolute value $<1$ (we recall that we are using the $t$-adic valuation). Having thus
infinitely many zeros in the domain of convergence and not being identically zero, $\phi$ is transcendental.

The series $\phi$ converges at $x=t$. Theorem \ref{theo1} implies that the formal series
$\phi(t)\in C((t))$ is transcendental over $\mathcal{K}$. We notice that the arguments of \ref{transcendencefunctionsvartheta} can be probably extended to investigate
the transcendence of the series $\phi$ associated to, say, $\vartheta=-(1+t/\theta)^{-1}$, case in which we do not necessarily have infinitely many zeros. The reason is that,
over $\KK((x))$, the $\FF_q$-linear Frobenius twist $F:a\mapsto a^q$ (for all $a$) splits as $$F=\tau\chi=\chi\tau$$ where $\tau$ is {\em Anderson's} $\FF_q((t))$-linear twist 
and $\chi$ is {\em Mahler's} $C((x))$-linear twist, and most of the arguments of \ref{transcendencefunctionsvartheta} can be generalised to this setting.


\subsubsection{Some $\tau$-difference equations in $\mathcal{K}[[x]]$}

The arguments of the previous subsection deal with formal series in $\mathcal{K}[[x]]=C(t)[[x]]$.
We have another important ring of formal series, also embedded in $\KK[[x]]$, which is $C(x)[[t]]$. Although
the arithmetic of values of these series seems to be not deducible from Theorem \ref{theo1},
we discuss here about some examples because solutions $f\in C(x)[[t]]$ of $\tau$-difference equations such as
\begin{equation}
\tau f=af+b,\quad a,b\in C(x)\label{abp}\end{equation} are often related to Anderson-Brownawell-Papanikolas
linear independence criterion in \cite{ABP} (see the corresponding contribution in this volume
and the refinement \cite{Chieh}).

With $\zeta_\theta$ a fixed $(q-1)$-th root of $-\theta$, the transcendental formal series
\begin{equation}
\boldsymbol{\Omega}(t):=\zeta_\theta^{-q}\prod_{i=1}^\infty\left(1-\frac{t}{\theta^{q^i}}\right)=\sum_{i=0}^\infty d_it^i\in K^\text{alg.}[[t]]
\label{eq:omega}
\end{equation}
is convergent for all $t\in C$, such that $\boldsymbol{\Omega}(\theta)\in\FF_q^\times\widetilde{\pi}^{-1}$ and satisfies
the functional equation
\begin{equation}
\boldsymbol{\Omega}(t^q)=(t^q-\theta^q)\boldsymbol{\Omega}(t)^q\label{omega}\end{equation} (see \cite{ABP}).
By a direct inspection it turns out that there is no finite extension of $K$ containing all the coefficients $d_i$ of the $t$-expansion of
$\boldsymbol{\Omega}$ (\footnote{They generate an infinite tower of Artin-Schreier extensions.}). Hence, there is no variant of Mahler's method which seems to apply to 
prove the transcendence of $\Omega$ at algebraic elements (and a suitable variant of Theorem \ref{corvajazannier} is not yet
available). 

The map $\tau:C((t))\rightarrow C((t))$ acts on in the following way: 
$$c=\sum_i c_it^i\mapsto \tau c:=\sum_i c_i^qt^i.$$ 

Let 
$s(t)$ be the formal series $\tau^{-1}\boldsymbol{\Omega}^{-1}\in C[[t]]$ (where $\tau^{-1}$ is the reciprocal map of $\tau$). After (\ref{omega})
this function is solution of the $\tau$-difference equation:
\begin{equation}\label{stheta}
(\tau s)(t)=(t-\theta)s(t),\end{equation} hence it is a solution of (\ref{abp}) with $a=t-\theta$ and $b=0$.
Transcendence of values of this kind of function does not seem to follow from Theorem \ref{theo1},
but can be obtained with \cite[Theorem 1.3.2]{ABP}.

More generally, let $\Lambda$ be an $A$-lattice of $C$ of rank $r$ and let
\begin{equation}
e_\Lambda(z)=z\prod_{\lambda\in\Lambda\setminus\{0\}}\left(1-\frac{z}{\lambda}\right)
\label{expo}\end{equation}
be its exponential function, in Weierstrass product form. The function $e_\Lambda$ is an entire, surjective, $\FF_q$-linear function. Let $\phi_\lambda:A\rightarrow\mathbf{End}_{\FF_q-\text{lin.}}(\GG_a(C))=C[\tau]$ (\footnote{Polynomial expressions in $\tau$ with the product satisfying $\tau c=c^q\tau$, for $c\in C$.}) be
the associated Drinfeld module. We have, for $a\in A$, $\phi_\Lambda(a)e_\Lambda(z)=e_\Lambda(az)$.

Let us extend $\phi_\Lambda$ over $\mathcal{K}$ by means of the endomorphism $\tau$ ($\tau t=t$). After having chosen an element $\omega\in \Lambda\setminus\{0\}$, define 
the function
\begin{equation}\label{eq:slambda}
s_{\Lambda,\omega}(t):=\sum_{n=0}^{\infty}e_\Lambda\left(\frac{\omega}{\theta^{n+1}}\right)t^n\in C[[t]],
\end{equation}
convergent for $|t|<q$. We have, for all $a\in A$,
\begin{equation}\label{stheta2}
\phi_\Lambda(a)s_{\Lambda,\omega}=\overline{a}s_{\Lambda,\omega},
\end{equation}
where, if $a=a(\theta)\in\FF_q[\theta]$, we have defined $\overline{a}:=a(t)\in\FF_q[t]$. This means that
$s_{\Lambda,\omega}$, as a formal series of $C[[t]]$, is an eigenfunction for all 
the $\FF_q((t))$-linear operators $\phi_\Lambda(a)$, with eigenvalue $\overline{a}$, for all $a\in A$.

 If $\Lambda=\widetilde{\pi}A$,
$\phi_\Lambda$ is Carlitz's module, and 
equation (\ref{stheta2}) implies the $\tau$-difference equation (\ref{stheta}).

\section{Algebraic independence \label{algindep}}

In \cite{Ma2}, Mahler proved his first result of algebraic independence obtained modifying and generalising 
the methods of his paper \cite{Ma1}. The result obtained involved $m$ formal series in several variables but we describe 
its consequences on the one variable theory only. Let $L$ be a number field embedded in $\CC$ and let $d>1$ be an integer.

\begin{Theoreme}[Mahler]
Given $m$ formal series $f_1,\ldots,f_m\in L[[x]]$,
satisfying functional equations
$$f_i(x^d)=a_if_i(x)+b_i(x),\quad 1\leq i\leq m$$
with $a_i\in L$, $b_i\in L(x)$ for all $i$, converging in the open unit disk.
If $\alpha$ is algebraic such that $0<|\alpha|<1$, then 
the transcendence degree over $L(x)$ of the field $L(x,f_1(x),\ldots,f_m(x))$ is equal to the 
transcendence degree over $\QQ$ of the field $\QQ(f_1(\alpha),\ldots,f_m(\alpha))$.
\end{Theoreme}

\subsection{Criteria of algebraic independence and applications.\label{criteria}}
Mahler's result remained nearly unobserved for several years. It came back to surface notably thanks to the work of Loxton and van der Poorten in the seventies, and 
then by Nishioka and several other authors. At the beginning, these authors developed criteria for algebraic independence tailored for application to algebraic independence of Mahler's values. Later, criteria for algebraic independence evolved in very general results, especially in the hands of Philippon. Here follows a particular case of a criterion of algebraic independence by Philippon The statement that follows merges the results \cite[Theorem 2]{Phi1} and \cite[Theorem 2.11]{Phi2} and uses
the data $\KK,\mathcal{L},|\cdot|,\mathcal{K},\mathcal{A},\ldots$ where $\KK$ is a complete algebraically closed field in two cases.
We examine only the cases in which $\KK$ is either $\CC$ or $C$, but it is likely that the principles of the criterion extend to 
several other fields, like that of Section \ref{foursteps}.

%

In the case $\KK=\CC$, $\mathcal{L}$ is a number field embedded in $\CC$, $|\cdot|$ is the usual absolute value, $h$ is the absolute logarithmic height
of projective points defined over $\QQ^{\text{alg.}}$, $\mathcal{K}$ denotes $\QQ$ and $\mathcal{A}$ denotes $\ZZ$. 

In the case $\KK=C$, $\mathcal{L}$ is a finite extension of $K=\FF_q(\theta)$ embedded in $C$, $|\cdot|$ is an absolute value associated to the $\theta^{-1}$-adic valuation,
$h$ is the absolute logarithmic height of projective points defined over $K^{\text{alg.}}$, $\mathcal{A}$ denotes the ring $A=\FF_q[\theta]$ and we write $\mathcal{K}=K$. 

In both cases,
if $P$ is a polynomial with coefficients in $\mathcal{L}$, we associate to it a projective point whose coordinates
are $1$ and its coefficients, and we write $h(P)$ for the logarithmic height of this point (which depends on $P$ up to permutation of the coefficients).

\begin{Theoreme}[Philippon]\label{theoremephilippon}
Let $(\alpha_1,\ldots,\alpha_m)$ be an element of $\KK^m$ and $k$ an integer with 
$1\leq k\leq m$. Let us suppose that
there exist three increasing functions $\ZZ_{\geq 1}\rightarrow\RR_{\geq1}$
\begin{eqnarray*}
\delta& & (\text{degree})\\
\sigma& & (\text{height})\\
\lambda& & (\text{magnitude})
\end{eqnarray*}
and five positive real numbers $c_2,c_3,c_\delta,c_\lambda,c_\sigma$
satisfying the following properties.
\begin{enumerate}
\item $\lim_{n\rightarrow\infty}\delta(n)=\infty$,
\item $\lim_{n\rightarrow\infty}\frac{s(n+1)}{s(n)}=c_s$ for $s=\delta,\sigma,\lambda$,
\item $\sigma(n)\geq\delta(n)$, for all $n$ big enough,
\item The sequence $n\mapsto\frac{\lambda(n)}{\delta(n)^{k+1}\sigma(n)}$ is ultimately increasing,
\item For all $n$ big enough, 
$$\lambda(n)^{k+1}>\sigma(n)\delta(n)^{k-1}(\lambda(n)^k+\delta(n)^k).$$
\end{enumerate}
Let us suppose that there exists a sequence of polynomials $(Q_n)_{n\geq 0}$
with $$Q_n\in \mathcal{L}[X_1,\ldots,X_m],$$ with $\deg_{X_i}Q_n\leq\delta(n)$ for all $i$ and $n$,
with $h(Q_n)\leq\sigma(n)$ for all $n$, with coefficients integral over $\mathcal{A}$, such that,
for all $n$ big enough,
$$-c_2\lambda(n)<\log|Q_n(\alpha_1,\ldots,\alpha_m)|<-c_3\lambda(n).$$ 
Then, the transcendence degree over $\mathcal{L}$ of the field $\mathcal{L}(\alpha_1,\ldots,\alpha_m)$
is $\geq k$.
\end{Theoreme}

Theorem \ref{theoremephilippon} can be applied to prove the following result. 
 
\begin{Theoreme}\label{theorem:denis} 
Let us assume that we are again in one of the cases above; $\KK=\CC$ or $\KK=C$, let $\mathcal{L}$ be as above.
Let $f_1,\ldots,f_m$ be formal series of $\mathcal{L}[[x]]$ an $d>1$ and integer, 
satisfying the following properties:
\begin{enumerate}
\item $f_1,\ldots,f_m$ converge for $|x|<1$,
\item $f_1,\ldots,f_m$ are algebraically independent over $\KK(x)$,
\item For all $i=1,\ldots,m$, there exist $a_i,b_i\in \mathcal{L}(x)$ such that
$$f_i(x^d)=a_i(x)f_i(x)+b_i(x),\quad i=1,\ldots,m.$$
\end{enumerate}
Let $\alpha\in\mathcal{L}$ be such that $0<|\alpha|<1$.
Then, for all $n$ big enough, $f_1(\alpha^{d^n}),\ldots,f_m(\alpha^{d^n})$ are algebraically independent
over $\mathcal{L}$.
\end{Theoreme}

This result is, for $\KK=\CC$, a corollary of Kubota's result \cite[Theorem p. 10]{Kubo}. For $\KK=C$, 
it is due to Denis \cite[Theorem 2]{Denis1}. See also \cite{Becker0, Denis2}.

\medskip

\noindent\emph{Sketch of proof of Theorem \ref{theorem:denis}
in the case $\KK=C$.} To simplify the exposition, we assume that $\mathcal{L}=K$. Let $N>0$ be an integer; 
there exists at least
one non-zero polynomial $P_N\in K[x,X_1,\ldots,X_n]$ (that we choose) of degree $\leq N$ in each indeterminate, 
such that the order $\nu(N)$ of vanishing at $x=0$ of the function $$F_N(x)=P_N(x,f_1(x),\ldots,f_m(x))$$
(not identically zero because of the hypothesis of algebraic independence of the functions
$f_i$ over $C(x)$), satisfies $\nu(N)\geq N^{m+1}$.

The choice of the parameter $N$ will be made later. If $c(x)\in A[x]$ is a non-zero polynomial such that
$ca_i,cb_i\in A[x]$ for $i=1,\ldots,m$, then we define, inductively, $R_0=P_N$ and
$$R_n=c(x)^NR_{n-1}(x^d,a_1X_1+b_1,\ldots,a_mX_m+b_m)\in A[x,X_1,\ldots,X_m].$$
Elementary inductive computations lead to the following estimates, holding for $n$ big enough depending on $N,f_1,\ldots,f_m$, where $c_4,c_5,c_6$ are integer constants 
depending on $f_1,\ldots,f_m$ only:
\begin{eqnarray*}
\deg_{X_i}(R_n)&\leq&N,\quad(i=1,\ldots,m),\\
\deg_x(R_n)&\leq&c_4d^nN,\\
\deg_\theta(R_n)&\leq&c_5d^nN,\\
h(R_n)&\leq&c_7+c_6d^nN,
\end{eqnarray*}
where we wrote $c_7(N)=h(R_0)$; it is a real number depending on $f_1,\ldots,f_m$ and $N$.

Since
$$R_n(x,f_1(x),\ldots,f_n(x))=(\prod_{i=0}^{n-1}c(x^{d^i})^N)R_0(x^{d^n},f_1(x^{d^n}),\ldots,f_m(x^{d^n})),$$
one verifies the existence of two constants $c_2>c_3>0$ such that
$$-c_2d^n\nu(N)\leq\deg_\theta(R_n(\alpha,f_1(\alpha),\ldots,f_m(\alpha)))\leq -c_3d^n\nu(N)$$
for all $n$ big enough, depending on $N,f_1,\ldots,f_m$ and $\alpha$.

Let us define:
$$Q_n(X_1,\ldots,X_n)=D^{c_1d^nN}R_n(\alpha,X_1,\ldots,X_m)\in A[\alpha][X_1,\ldots,X_m],$$
where $D\in A\setminus\{0\}$ is such that $D\alpha$ is integral over $K$.
The estimate above implies at once that, for $n$ big enough:
\begin{eqnarray*}
\deg_{X_i}Q_n&\leq &N,\\
h(Q_n)&\leq &c_8(N)+c_9d^nN,
\end{eqnarray*}
where $c_8(N)$ is a constant depending on $f_1,\ldots,f_m,N$, and $\alpha$  (it can be computed with an explicit dependence on $c_7(N)$)
and $c_9$ is a constant depending on $f_1,\ldots,f_m$, and $\alpha$ but not on $N$.

Finally, Theorem \ref{theoremephilippon} applies with the choices:
\begin{eqnarray*}
\alpha_i&=&f_i(\alpha),\quad i=1,\ldots,m\\
k&=&m\\
\lambda(n)&=& d^n\nu(N)\\
\delta(n)&=&N\\
\sigma(n)&=&c_9d^n\nu(N),
\end{eqnarray*}
provided that we choose $N$ large enough depending on the constants $c_1,\ldots$ introduced so far. Then, one chooses $n$ big enough (depending on the good choice of $N$).\CVD

\subsubsection{An example with complex numbers.\label{examplecomplex}}

Theorem \ref{theorem:denis} furnishes algebraic independence of Mahler's values if we are able to check algebraic independence of Mahler's functions but it does not say anything on the latter problem; this is not an easy task in general. With 
the following example, we would like to
sensitise the reader to this problem which, the more we get involved in the subtleties of Mahler's method,
the more it takes a preponderant place.

In the case $\KK=\CC$, we consider the formal series in $\ZZ[[x]]$:
$$L_0=\prod_{i=0}^\infty(1-x^{2^i})^{-1},\quad L_r=\sum_{i=0}^\infty\frac{x^{2^ir}}{\prod_{j=0}^{i-1}(1-x^{2^j})},\quad r\geq 1$$
converging, on the open unit disk $|x|<1$, to functions satisfying~:
$$L_0(x^2)=(1-x)L_0(x),\quad L_r(x^2)=(1-x)(L_r(x)-x^r),\quad r\geq 1.$$ 

By a result of Kubota \cite[Theorem 2]{Kubo}, (see also Tšpfer,
\cite[Lemma 6]{Toepfer}), if $(L_i)_{i\in\mathcal{I}}$
(with $\mathcal{I}\subset\ZZ$) are algebraically dependent, then they also
are
$\CC$-linearly dependent modulo $\CC(x)$ in the following sense. There
exist complex numbers $(c_i)_{i\in\mathcal{I}}$
not all zero, such that:
$$\sum_{i=1}^mc_iL_i(x)=f(x)$$ with $f(x)\in\CC(x)$.
P. Bundschuh pointed out
that $L_0,L_2,L_4,\ldots$ are algebraically independent. To obtain this
property, he studied the behaviour of these functions near the the unit
circle. For a long time the author was convinced of the algebraic
independence of the functions $L_0,L_1,\ldots$ until very recently, when
T. Tanaka and H. Kaneko
exhibited non-trivial linear relations modulo $\CC(x)$ involving
$L_0,\ldots,L_s$ for all $s\geq 1$, some of which looking very simple, such as the relation:
$$L_0(x)-2L_1(x)=1-x.$$

\subsection{Measuring algebraic independence}

Beyond transcendence and algebraic independence, the next step in the study of the arithmetic of
Mahler's numbers is that of quantitative results such as measures of algebraic independence. Very often, 
such estimates are not mere technical refinements of well known results but deep information on the 
diophantine behaviour of classical constants; everyone knows the important impact that Baker's theory on quantitative minorations
of linear forms in logarithms on algebraic groups had in arithmetic geometry.

Rather sharp estimates are known for complex Mahler's values.
We quote here a
result of Nishioka \cite{Ni2, Ni} and \cite[Chapter 12]{Nesterenko:Introduction12} (it has been generalised by Philippon: \cite[Theorem 6]{Phi3}).

\begin{Theoreme}[Nishioka]\label{measurenishioka} Let us assume that, in the notations previously introduced, $\KK=\CC$. Let $L$ 
be a number field embedded in $\CC$.
Let $f_1,\ldots,f_m$ be formal series of $L[[x]]$, let us write $\underline{f}\in\mathbf{Mat}_{n\times 1}(L[[x]])$ for the column matrix
whose entries are the $f_i$'s. Let
$\mathcal{A}\in\mathbf{Mat}_{n\times n}(L(x)), 
\underline{b}\in\mathbf{Mat}_{n\times 1}(L(x))$ be matrices.
Let us assume that:
\begin{enumerate}
\item $f_1,\ldots,f_m$ are algebraically independent over $\CC(x)$,
\item For all $i$, the formal series $f_i(x)$ converges for $x$ complex such that $|x|<1$,
\item $\underline{f}(x^d)=\mathcal{A}(x)\cdot\underline{f}(x)+\underline{b}(x)$.
\end{enumerate}

Let $\alpha\in L$ be such that $0<|\alpha|<1$, not a zero or a pole of $\mathcal{A}$ and not a pole of $\underline{b}$.
Then, there exists a constant $c_1>0$ effectively computable
depending on $\alpha,\underline{f}$, with the following property.

For any $H,N\geq 1$ integers and any non-zero polynomial
$P\in\ZZ[X_1,\ldots,X_m]$ whose partial degrees in every indeterminate do not exceed $N$ and whose coefficients are not greater that $H$
in absolute value, the number $P(f_1(\alpha),\ldots,f_m(\alpha))$ is non-zero and the inequality below holds:
\begin{equation}\label{estimationnishioka}
\log|P(f_1(\alpha),\ldots,f_m(\alpha))|\geq -c_1N^m(\log H+N^{m+2}).\end{equation}
\end{Theoreme}

We sketch how Theorem \ref{measurenishioka} implies the algebraic independence of $f_1(\alpha),\ldots,f_{m}(\alpha)$ and
$g(\alpha)$ with $f_1,\ldots,f_m,g\in \QQ[[x]]$ algebraically independent over $\QQ[[x]]$ satisfying linear functional equations
as in Theorem \ref{theorem:denis} and $g$ satisfying $$g(x^d)=a(x)g(x)+b(x),$$ with $a,b\in \QQ(x)$, $\mathcal{A},\underline{b}$ with rational coefficients, and $\alpha$ not a pole of all these rational
functions. 
Of course, this is a simple corollary of Theorem \ref{measurenishioka}. However, we believe that the proof is instructive;
it follows closely Philippon's ideas in \cite{Phi3}. The result is reached because the estimates of Theorem \ref{measurenishioka} are precise enough.
In particular, the separation of the quantities $H$ and $N$ in (\ref{estimationnishioka}) is crucial.

\subsubsection{Algebraic independence from measures of algebraic independence.\label{frommeasuretoindep}}

For the purpose indicated at the end of the last subsection, we assume that $\alpha\in\QQ^\times$. This hypothesis in not strictly necessary and is assumed only to 
simplify the exposition of the proof; by the way, the reader will remark that several other hypotheses we assume are avoidable.

\medskip

\noindent{\em Step (AP)}. For all $N\geq 1$, we choose a non-zero polynomial $P_N\in\ZZ[x,X_1,\ldots,X_m,Y]$ of partial degrees $\leq N$ in each indeterminate,
such that, writing $$F_N(x):=P_N(x,f_1(x),\ldots,f_m(x),g(x))=c_{\nu(N)}x^{\nu(N)}+\cdots\in\QQ[[x]],\quad c_{\nu(N)}\neq0,$$ we have
$\nu(N)\geq N^{m+2}$ (we have already justified why such a kind of polynomial exists).

Just as in the proof of Theorem \ref{theorem:denis} we construct, for each $N\geq 1$, a sequence of polynomials $(P_{N,k})_{k\geq 0}$
in $\ZZ[x,X_1,\ldots,X_m,Y]$
recursively in the following way:
\begin{eqnarray*}
P_{N,0}&:=&P_N,\\
P_{N,k}&:=&c(x)^NP_{N,k-1}(x^d,a_1(x)X_1+b_1(x),\ldots,a_m(x)X_m+b_m(x),a(x)Y+b(x)),
\end{eqnarray*}
where $c(x)\in\ZZ[x]\setminus\{0\}$ is chosen so that $ca_i,cb_j,ca,cb$ belong to $\ZZ[x]$.
The following estimates are easily obtained:
\begin{eqnarray*}
\deg_ZP_{N,k}&\leq &N,\quad \text{ for }Z=X_1,\ldots,X_m,Y,\\
\deg_xP_{N,k}&\leq & c_2d^kN,\\
h(P_{N,k})&\leq &c_3(N)+c_4d^kN,\\
\end{eqnarray*}
where $c_2,c_4$ are positive real numbers effectively computable depending on $\alpha,\underline{f}$
and $g$, and $c_3(N)>0$ depends on these data as well as on $N$ (it depends on the choice of the polynomials $P_N$).

Let us assume by contradiction that $g(\alpha)$ is algebraic over the field 
$$\mathcal{F}:=\QQ(f_1(\alpha),\ldots,f_m(\alpha)),$$
of transcendence degree $m$ over $\QQ$.
We observe that, after the identity principle of analytic functions we have, for $k$ big enough depending on $\alpha,\underline{f},g$ and $N$:
\begin{equation}\label{identityprincipleagain}
P_{N,k}(\alpha,f_1(\alpha),\ldots,f_m(\alpha),g(\alpha))\in\mathcal{F}^\times.\end{equation}
Let $\widetilde{Q}\in\mathcal{F}[X]\setminus\{0\}$ be the minimal polynomial of $g(\alpha)$, algebraic over $\mathcal{F}$.
We can write $\widetilde{Q}=a_0+a_1X+\cdots+a_{r-1}X^{r-1}+X^r$ with the $a_i$'s in $\mathcal{F}$. Multiplying
by a common denominator, we obtain a non-zero polynomial $Q\in\ZZ[X_1,\ldots,X_m,Y]$ such that 
$Q(f_1(\alpha),\ldots,f_m(\alpha),g(\alpha))=0$, with the property that the polynomial $Q^*=Q(f_1(\alpha),\ldots,f_m(\alpha),Y)\in\mathcal{F}[Y]$
is irreducible.

\medskip

\noindent{\em Step (NV)}. Let us denote by $\Delta_k$ the resultant $\mathbf{Res}_Y(P_{N,k},Q)\in\ZZ[x,X_1,\ldots,X_m]$.
If $\delta_k:=\Delta_k(\alpha,f_1(\alpha),\ldots,f_m(\alpha))\in\mathcal{F}$ vanishes for a certain $k$,
then $Q^*$ and 
$$P^*_{N,k}:=P_{N,k}(\alpha,f_1(\alpha),\ldots,f_m(\alpha),Y)\in\mathcal{F}[Y]$$ have a
common zero. Since $Q^*$ is irreducible, we have that $Q^*$ divides $P^*_{N,k}$ in $\mathcal{F}[Y]$
and $$P_{N,k}(\alpha,f_1(\alpha),\ldots,f_m(\alpha),g(\alpha))=0;$$ this cannot happen for $k$ big enough
by the identity principle of analytic functions (\ref{identityprincipleagain}) so that we can assume that
for $k$ big enough, $\delta_k\neq 0$, ensuring that $\Delta_k$ is not identically zero; the estimates  of
the height and the degree of $\Delta_k$ quoted below are simple exercises and we do not give details of their proofs:
\begin{eqnarray*}
\deg_Z\Delta_k&\leq &c_5N,\quad Z=X_1,\ldots,X_m,\\
\deg_x\Delta_k&\leq &c_6d^kN,\\
h(\Delta_k)&\leq &c_7(N)+c_8d^kN,
\end{eqnarray*}
where $c_5,c_6,c_8$ are positive numbers effectively computable depending on $\alpha,\underline{f}$ and $g$, while
the constant $c_7(N)$ depends on these data and on $N$. 

Let $D$ be a non-zero positive integer such that $D\alpha\in\ZZ$. Then, writing $$\Delta^*_k:=D^{Nd^k}\Delta_k(\alpha,X_1,\ldots,X_m),$$
we have $\Delta_k^*\in\ZZ[X_1,\ldots,X_m]\setminus\{0\}$ and
\begin{eqnarray*}
\deg_{X_i}\Delta^*_k&\leq& c_5N,\\
h(\Delta_k^*)&\leq &c_7(N)+c_8d^kN.
\end{eqnarray*}

\noindent{\em Step (LB)}.
By Nishioka's Theorem \ref{measurenishioka}, we have the inequality (for $k$ big enough):
\begin{equation}\label{minorationnishioka}
\log|\Delta_k^*(f_1(\alpha),\ldots,f_m(\alpha))|\geq c_1N^m(d^kN+c_9(N)),
\end{equation}
where $c_9(N)$ is a constant depending on $N$.

To finish our proof, we need to find an upper bound contradictory with (\ref{minorationnishioka}); it will be obtained 
by analytic estimates as usual.

\medskip

\noindent{\em Step (UB)}. Looking at the proof of Lemma 5.3.1 of \cite{Wa2}, and using in particular inequality (1.2.7) of
loc. cit., we verify the existence of constants $c_{12},c_{15}$ depending on $\alpha,\underline{f},g$,
$c_{14}(\alpha,\epsilon)$ depending on $\alpha$ and $\epsilon$, and $c_{13}(N)$ depending on $\alpha,\underline{f},g$ and $N$, such that:
\begin{eqnarray}
\log|\Delta_k^*(f_1(\alpha),\ldots,f_m(\alpha))|&\leq & \log(N+c_{10})+c_{11}h(P_{N,k}^*)+(N+1)h(Q)+\nonumber\\
& &
\log|c_{\nu_N}|+\nu(N)d^k\log|\alpha|+\epsilon\nonumber\\
&\leq &c_{12}(\log N+d^kN)+c_{13}(N)-c_{14}(\alpha,\epsilon)d^k\nu(N)\nonumber\\
&\leq & c_{15}d^kN+c_{13}(N)-c_{14}(\alpha,\epsilon)d^kN^{m+2}.\label{upperboundnishioka}
\end{eqnarray}
Finally, it is easy to choose $N$ big enough, depending on $c_1,c_{15},c_{14}$ but {\em not} on $c_9,c_{13}$
so that, for $k$ big enough, the estimates (\ref{minorationnishioka}) and (\ref{upperboundnishioka}) are not compatible: this is due to the particular shape of
(\ref{estimationnishioka}), with the linear dependence in $\log H$. 

\subsubsection{Further remarks, comparisons with Nesterenko's Theorem.}

In the sketch of proof of the previous subsection, the reader probably observed a kind of induction structure; a measure of algebraic independence for $m$ numbers delivers
algebraic independence for $m+1$ numbers. The question is then natural: {\em is it possible to obtain a measure of algebraic independence for $m+1$ numbers
allowing continue the process and consider $m+2$ numbers?}

In fact {\em yes}, there always is an inductive structure of proof, but {\em no}, it is not just a measure for $m$ numbers which alone
implies a measure for $m+1$ numbers.
Things are more difficult than they look at first sight
and the inductive process one has to follow concerns other parameters as well. For instance, the reader can verify that  
it is unclear how to generalise the arguments of \ref{frommeasuretoindep} and work directly with a polynomial $Q$ which has a very small value at $\omega=(f_1(\alpha),\ldots,f_m(\alpha),g(\alpha))$.

Algebraic independence theory usually appeals to {\em transfer techniques}, as an alternative to direct estimates at $\omega$. 
A detour on a theorem of Nesterenko might be useful to understand what is going on so our discussion now temporarily leaves
Mahler's values, that will be reconsidered in a little while.

Precise multiplicity estimates in differential rings generated by Eisenstein's series obtained by Nesterenko,
the criteria for algebraic independence by Philippon already mentioned in this paper and a trick 
of the stephanese team (cf. \ref{hauptmodul})
allowed Nesterenko, 
in 1996, to prove the following theorem
(see \cite{Nest,Nesterenko:Introduction3}):

 \begin{Theoreme}[Nesterenko]\label{nesterenkotheorem}
Let $E_2,E_4,E_6$ the classical {\em Eisenstein's series of wei\-ghts $2,4,6$} respectively, normalised so that 
$\lim_{\Im(z)\rightarrow\infty}E_{2i}(z)=1$ (for $\Re(z)$ bounded), let $z$ be a complex number of strictly positive imaginary part.
Then, three of the four complex numbers $e^{2\pi\mathrm{i}z},E_2,E_4,E_6$ are algebraically independent. 
\end{Theoreme}

Although Eisenstein's series are not directly related to Mahler's functions there is a hidden link and the ideas introduced to 
prove Theorem \ref{nesterenkotheorem} influenced Nishioka in her proof of Theorem \ref{estimationnishioka} as well as 
other results by Philippon that we will mention below. This is why we cannot keep silent on this aspect.

First of all, we recall that in \cite{Phi3}, Philippon
showed how to deduce the algebraic independence of $\pi,e^{\pi},\Gamma(1/4)$ ({\footnote{Philippon's result is in fact more general than the algebraic independence of these three numbers, but less general than Nesterenko's theorem \ref{nesterenkotheorem}, although it historically followed it. Our arguments in \ref{frommeasuretoindep}
are strongly influenced by it. In \cite{Phi3}, Philippon proposes alternative, simpler proofs for Nesterenko's theorem.}), a
well known corollary of Nesterenko's Theorem \ref{nesterenkotheorem}, from a measure of algebraic independence of $\pi,\Gamma(1/4)$ by Philibert
in \cite{Philibert}.
This implication was possible because Philibert's result was sharp enough.
It is however virtually impossible to deduce Theorem \ref{nesterenkotheorem} or the quantitative result in \cite{Nest} which can also be deduced from corollary of \cite[Theorem 3]{Phi3}
(\footnote{A result asserting that, for a polynomial $P\in\ZZ[X_0,\ldots,X_4]\setminus\{0\}$ 
whose partial degrees in every indeterminate do not exceed $N>0$ and whose coefficients are not greater that $H>0$
in absolute value, then $\log|P(e^{\pi},\pi,\Gamma(1/4))|\geq -c_1(\epsilon)(N+\log H)^{4+\epsilon}$, where $c_1$ is an absolute constant depending on $\epsilon$ only.
The dependence in $\epsilon$ is completely explicit.})
just by using Philibert's result.

In the proof of Theorem \ref{estimationnishioka}, Nishioka proves (just as Nesterenko does in \cite{Nest}) a more general
measure of the smallness of the values that a polynomial with rational integer coefficients assumes at $(1,f_1(\alpha),\ldots,f_m(\alpha),g(\alpha))$,
restricting the choice of that polynomial  in a given unmixed homogeneous ideal $I$. The proof of such a kind of result (cf. \cite[Lemma 2.3]{Nesterenko:Introduction12})
involves induction on the dimension of $I$.

Assuming the existence of an ideal $I$ with minimal dimension ``very small" at $\omega$, it is possible, looking at its reduced primary decomposition, to 
concentrate our attention to $I=\mathfrak{p}$ prime. The ``closest point principle" of \cite[p. 89]{Nesterenko:Introduction6} allows to show the existence,
in the projective variety $V$ associated to $\mathfrak{p}$, of a point $\beta$ which is at a very short distance from $\omega$ (see also \cite[Proposition 1.5]{Nesterenko:Introduction3}). This shows that in this problem, to measure a polynomial or an ideal at $\omega$ it is more advantageous to do it at $\beta$;
indeed, all the polynomials of $\mathfrak{p}$ vanish at $\beta$.

At this point, it remains to construct an unmixed homogeneous ideal $J$ of dimension $\dim\mathfrak{p}-1$, contradicting our assumptions.
To do so, it is necessary to proceed as we did in \ref{frommeasuretoindep} to construct $P_{N,k}$ etc., with the important difference that now, 
all the estimate depend very much on the choice, that must then use Siegel's Lemma. Another important tool that has to be used is 
a  {\em multiplicity estimate}, that belongs to step (NV), proved by Nishioka \cite[Theorem 4.3]{Ni}, that we reproduce here.

\begin{Theoreme}[Nishioka]\label{nesterenkonishioka} Let $f_1,\ldots,f_m$ be satisfying the hypotheses of Theorem \ref{measurenishioka}, so
that for all $P\in\CC[x,X_1,\ldots,X_m]\setminus\{0\}$, the function $F(x):=P(x,f_1(x),\ldots,f_m(x))$ 
has the expansion
$$F(x)=c_{\nu}x^{\nu}+\cdots,$$ with $c_{\nu}\neq0$.
There exists a constant $c_1>0$, depending on $f_1,\ldots,f_m$ only, with the following property.
If $P$ is as above and $N_1:=\max\{1,\deg_xP\}$ and $N_2:=\max\{1,\deg_{X_1}P,\ldots,\deg_{X_m}P\}$,
then
$$\nu\leq c_1N_1N_2^m.$$
\end{Theoreme} 
This result is very similar to Nesterenko's multiplicity estimate \cite[Chapter 10, Theorem 1.1]{Nesterenko:Introduction3} and again,
its proof essentially follows Nesterenko's ideas.

The ideal $J$ previously mentioned  is defined as the ideal generated by $\mathfrak{p}$ and a polynomial obtained from $P_{N,k}$ by homogenisation, substitution $x=\alpha$,
and a good choice of $N,k$ taking into account the magnitude of the coefficients of the series $f_i,g$. Indeed, one proves that such a polynomial cannot belong to $\mathfrak{p}$. The closest point principle is necessary in this kind of proof. 

The arguments of the above discussion can be modified to obtain the analog of Theorem \ref{nesterenkotheorem}
for values of Mahler's functions at general complex numbers, obtained by Philippon (cf. \cite[Theorem 4]{Phi3}).
Here, $L$ is again a number field embedded in $\CC$ and $d>1$ is an integer.

\begin{Theoreme}[Philippon]\label{theophili} Under the same hypotheses and notations of Theorem \ref{measurenishioka},
if $\alpha$ is a complex number with $0<|\alpha|<1$, then,
for $n$ big enough, the complex numbers $\alpha,f_1(\alpha^{d^n}),\ldots,f_m(\alpha^{d^n})\in\CC$ generate a subfield of $\CC$ of transcendence degree $\geq m$. 
\end{Theoreme}
This result is a corollary of a more general quantitative result \cite[Theorem 6]{Phi3} which follows from
Philippon's criterion for measures of algebraic independence (loc. cit. p. 5). 

\subsubsection{Commentaries on the case of positive characteristic}

Similar, although simpler arguments are in fact commonly used to obtain 
measures of transcendence. Several authors deduce them from measures of linear algebraic approximation; see for example
Amou, Galochkin and Miller \cite{Amou, Gal,Mil} (\footnote{Some results hold for series which satisfy functional equations  which are not necessarily linear.}).
These results often imply that Mahler's values are Mahler's {\em $S$-numbers}. 

In positive characteristic, it is well known that  separability difficulties occur preventing to deduce good measures of transcendence from 
measures of linear algebraic approximation (\footnote{As first remarked by Lang, for complex numbers, there is equivalence between measures of transcendence and measures of linear algebraic approximation in the sense that, from a measure of linear algebraic approximation one can get a measure of transcendence and then again,
a measure of linear algebraic approximation which is essentially that of the beginning, with a controllable degradation of the constants.}).
In \cite{denis3}, Denis proves the following result, where $d>1$ is an integer.

\begin{Theoreme}[Denis]\label{nishioka1} Let us consider a finite extension $\mathcal{L}$ of $K=\FF_q(\theta)$, $f\in \mathcal{L}[[x]]$, $\alpha\in \mathcal{L}$ such that $0<|\alpha|<1$. Let us assume that
$f$ is transcendental over $C(x)$, convergent for $x\in \mathcal{C}$ such that $|x|< 1$, and satisfying the linear functional equation
$$f(x^q)=a(x)f(x)+b(x),$$ with $a,b\in \mathcal{L}(x)$.

For all $n$ big enough, we have the following property. Let 
$\beta=\alpha^{d^n}$. Then, there exists an effectively computable constant $c_1>0$
depending on $\beta,f$ only, such that,
given any non-constant polynomial $P\in A[X]$,
\begin{equation}\label{miller2}
\log|P(f(\beta))|\geq -\deg_X(P)^4(\deg_X(P)+\deg_\theta(P)).
\end{equation}
\end{Theoreme}
This result yields completely explicit measures of transcendence of $\widetilde{\pi}$, 
of Carlitz's logarithms of elements of $K$, and of certain Carlitz-Goss's zeta values (see Section \ref{carlitz} 
for definitions).

To prove Theorem \ref{nishioka1}, Denis uses the following multiplicity estimate.
\begin{Theoreme}\label{multiplicityestdenis}
Let $\mathcal{K}$ be any (commutative) field and $f\in \mathcal{K}((x))$ be transcendental satisfying the 
functional equation
$f(x^d)=R(x,f(x))$ with $R\in \mathcal{K}(X,Y)$ and $h_Y(R)<d$.
Then, if $P$ is a polynomial in $\mathcal{K}[X,Y]\setminus\{0\}$ such that
$\deg_XP\leq N$ and $\deg_YP\leq M$, we have
$$\text{ord}_{x=0}P(x,f(x))\leq N(2Md+Nh_X(Q)).$$ 
\end{Theoreme}
It would be interesting to generalise such a multiplicity estimate for several
algebraically independent formal series and obtain a variant of Tšpfer's \cite[Theorem 1]{Toepfer2} (see \cite{Pellarin1, Pellarin2, Pellarin3} to check the difficulty 
involved in the research of an analogue of Nesterenko's multiplicity estimate
for Drinfeld quasi-modular forms). This could be helpful to obtain analogues of Theorem \ref{measurenishioka} for Mahler's values 
in fields of positive characteristic.

\subsection{Algebraic independence of Carlitz's logarithms\label{carlitz}}

For the rest of this chapter, we will give some application of Mahler's method and of Anderson-Brownawell-Papanikolas method to algebraic independence of Carlitz's logarithms of algebraic elements of $C$
and of some special values of Carlitz-Goss zeta function at rational integers.
Results of this part are not original since they are all contained in the papers \cite{Papanikolas1} by Papanikolas
and \cite{ChangYu} by Chieh-Yu Chang and Jing Yu.
But the methods we use here are slightly different and self-contained. 

Both proofs of the main results in \cite{Papanikolas1, ChangYu} make use of a general statement \cite[Theorem 5.2.2]{Papanikolas1} which can be considered as a
variant of Grothendieck period conjecture for a certain generalisation of Anderson's $t$-motives, also due Papanikolas. To apply this result,  the computation of motivic Galois groups associated
to certain $t$-motives is required.

Particular cases of these results are also contained in Denis work \cite{denis},
where he applies Mahler's method and without appealing to Galois theory.
Hence, we follow the ideas of the example in \ref{examplecomplex}
and the main worry here is to develop analogous proofs in the Drinfeldian framework.
In \ref{examplecomplex} the explicit computation of the transcendence degree 
of the field generated by $L_0,L_1,\ldots$ was pointed out as a problem.
But we have already remarked there, that if $L_0,L_1,\ldots$ are algebraic dependent,
then they also are $\CC$-linearly dependent modulo $\CC(x)$.
This property, consequence of a result by Kubota, is easy to obtain because the matrix 
of the linear difference system of equations satisfied by the $L_i$'s has the matrix of
its associated homogeneous system which is {\em diagonal}.

For $\Lambda=\widetilde{\pi}A$ with $\widetilde{\pi}$ as in (\ref{prodottopi}),
the exponential function $e_{\carlitz}:=e_\Lambda$ (\ref{expo})
can be explicitly written as follows:
$$e_{\carlitz}(z)=\sum_{i\geq 0}\frac{z^{q^i}}{[i][i-1]^q\cdots[1]^{q^{i-1}}},
$$
where $[i]:=\theta^{q^i}-\theta$ ($i\geq 1)$.
This series converges uniformly on every open ball with center in $0$
to an $\FF_q$-linear surjective function $e_\carlitz:C\rightarrow C$.
The formal series $\log_{\carlitz}$, reciprocal of $e_\carlitz$  in $0$, converges for $|z|<q^{q/(q-1)}=|\widetilde{\pi}|$.
Its series expansion can be computed explicitly:
$$\log_{\carlitz}(z)=\sum_{i\geq 0}
\frac{(-1)^iz^{q^i}}{[i][i-1]\cdots[1]}.$$

The first Theorem we shall prove in a simpler way is the following (cf. \cite[Theorem 1.2.6]{Papanikolas1}):

\begin{Theoreme}[Papanikolas]\label{papa1}
Let $\ell_1, \dots, \ell_m \in C$ be such that $e_{\carlitz}(\ell_i)
\in K^{\text{alg.}}$ ($i=1, \dots, m$).  If $\ell_1, \dots, \ell_m$
are linearly independent over $K$, then they also are algebraically independent over $K$.\end{Theoreme}

\subsubsection{Carlitz-Goss polylogarithms and zeta functions.}
Let us write $A_+=\{a\in A, \; a \textrm{ monic}\}$.
In \cite{GossZeta}, Goss introduced a function $\zeta$,
defined over $C\times\ZZ_p$ with values in $C$,
such that for $n\geq 1$ integer, $$\zeta(\theta^n,n)=\sum_{a\in A_+}\frac{1}{a^n}\in K_\infty.$$
In the following, we will write $\zeta(n)$ for $\zeta(\theta^n,n)$.
For $n\in\NN$, let us also write
$\Gamma(n):=\prod_{i=0}^sD_i^{n_i}\in K$, $n_0+n_1q+\cdots+n_sq^s$ being
the expansion of $n-1$ in base $q$
and $D_i$ being the polynomial $[i][i-1]^q\cdots[1]^{q^{i-1}}$.
It can be proved that $z/e_{\carlitz}(z)=\sum_{n=0}^\infty
B_n\frac{z^n}{\Gamma(n+1)}$ for certain $B_{n}\in K$. The so-called Bernoulli-Carlitz
relations can be obtained by a computation involving the logarithmic derivative of $e_\carlitz(z)$:
for all $m\geq 1$,
\begin{equation}
\frac{\zeta(m(q-1))}{\widetilde{\pi}^{m(q-1)}}=\frac{B_m}{\Gamma(m(q-1)+1)}\in K.\label{eq:zeta_bernoulli}\end{equation}
In particular, one sees that
\begin{eqnarray*}
\widetilde{\pi}^{q-1}&=&(\theta^q-\theta)\zeta(q-1)\in K_\infty.\label{eq1}
\end{eqnarray*}
We also have the obvious relations:
\begin{equation}\label{eq:relations_triviales}
\zeta(mp^k)=\zeta(m)^{p^k},\quad m,k\geq 1.
\end{equation}
The second theorem we are going to prove directly is:

\begin{Theoreme}[Chang, Yu]\label{changyu}
The algebraic dependence relations over $K$ between the numbers $$\zeta(1),\zeta(2),\ldots$$
are generated by Bernoulli-Carlitz's relations
(\ref{eq:zeta_bernoulli}) and the relations (\ref{eq:relations_triviales}).
\end{Theoreme}

\subsubsection{Two propositions.}

In this subsection we develop an analogue of \cite[Lemma 6]{Toepfer}, for the same purpose
we needed it in  \ref{examplecomplex}. 

We consider here a perfect field $U$ of characteristic $p>0$ containing $\FF_q$ 
and a $\FF_q$-auto\-morphism $\tau:U\rightarrow U$. Let $U_0$ be the subfield of constants of $\tau$,
namely, the subset of $U$ whose elements $s$ are such that $\tau s=s$.

For example, we can consider $U=\bigcup_{n\geq 0}C(x^{1/p^n})$ with $\tau$ defined 
as the identity over $C$, with $\tau x=x^q$. Another choice is to consider $U=\bigcup_{n\geq 0}C(t^{1/p^n})$,
with $\tau$ defined by $\tau c=c^{1/q}$ for all $c\in C$ and $\tau t=t$. In the first example we have $U_0=C$ while in the second, $U_0=\bigcup_{n\geq 0}\FF_q((t^{1/p^n}))$.

More generally, after \ref{foursteps}, we can take either
$U=\bigcup_{n\geq 0}\KK(x^{1/p^n})$ or the field $\bigcup_{n,m\geq 0}C(x^{1/p^n},t^{1/p^m})$ (which is contained
in the previous field)
with the corresponding automorphism $\tau$ (these settings will essentially include the two examples above).
In the first case, we have $U_0=\FF_q\langle\langle t\rangle\rangle$, and in the second case,
we have $U_0=\bigcup_{n\geq 0}\FF_q((t^{1/p^n}))$.

Let us also consider the ring ${\cal R}=U[X_1,\ldots,X_N]$ and write, 
for a polynomial $P=\sum_{\underline{\lambda}}c_{\underline{\lambda}}\underline{X}^{\underline{\lambda}}\in{\cal R}$, $P^\tau$ as the polynomial $\sum_{\underline{\lambda}}(\tau c_{\underline{\lambda}})\underline{X}^{\underline{\lambda}}$.
Let $D_1,\ldots,D_N$ be elements of $U^\times$, $B_1,\ldots,B_N$ be elements of $U$ and, 
for a polynomial $P\in{\cal R}$, let us write 
$$\widetilde{P}=P^\tau(D_1X_1+B_1,\ldots,D_NX_N+B_N).$$ 

We prove the following two Propositions, which provide together  the analogue in positive characteristic
of Kubota \cite[Theorem 2]{Kubo}.

\begin{Proposition}\label{propo1}
Let $P\in{\cal R}$ be a non-constant polynomial such that $\widetilde{P}/P\in{\cal R}$.
Then, there exists a polynomial $G\in{\cal R}$ of the form $G=\sum_ic_iX_i+B^p$
such that $\widetilde{G}/G\in{\cal  R}$,
where $c_1,\ldots,c_N\in U$ are not all vanishing and $B\in{\cal R}$. If $W$ is the subfield
generated by $\FF_q$ and the coefficients of $P$, then there exists $M\geq 1$
such that for each coefficient $c$ of $G$, $c^{p^M}\in W$.
\end{Proposition}

\noindent\emph{Proof.} 
If $P\in{\cal R}$ is such that $\widetilde{P}=QP$ for $Q\in {\cal R}$ one sees, comparing the degrees of $\widetilde{P}$
and $P$, that $Q\in U$ and if $P$ is non-zero, $Q\not=0$. The subset of ${\cal R}$ of these
polynomials is a semigroup ${\cal S}$ 
containing $U$.
If $P\in{\cal S}$ satisfies $\widetilde{P}=QP$, then $F:=\partial P/\partial X_i$
belongs to ${\cal  S}$ since $\widetilde{F}=D_i^{-1}QF$. Similarly, if 
$P=F^p\in{\cal S}$ with $F\in{\cal S}$
then $F\in{\cal S}$ as one sees easily that in this case, $\widetilde{F}=Q^{1/p}F$. 

By hypothesis, ${\cal S}$ contains a non-constant polynomial $P$.
We now show that the polynomial $G\in{\cal S}$ as in the Proposition can be constructed 
by iterated applications
of partial derivatives $\partial_1=\partial/\partial X_1,\ldots,\partial_N=\partial/\partial X_N$ and $p$-root extrations starting from $P$. 

Let $P$ be as in the hypotheses. We can assume that $P$ is not a $p$-th power.
We can write:
$$P=\sum_{\underline{\lambda}=(\lambda_1,\ldots,\lambda_N)\in\{0,\ldots,p-1\}^N}c_{\underline{\lambda}}\underline{X}^{\underline{\lambda}},\quad c_{\underline{\lambda}}\in{\cal R}^p.$$
Let $M:=\max\{\lambda_1+\cdots+\lambda_N,\; c_{\underline{\lambda}}\not=0\}$. We can write
$P=P_1+P_2$ with $$P_1:=\sum_{\lambda_1+\cdots+\lambda_N=M}c_{\underline{\lambda}}\underline{X}^{\underline{\lambda}}.$$
There exists $(\beta_1,\ldots,\beta_N)\in\{0,\ldots,p-1\}^N$ 
with $\beta_1+\cdots+\beta_N=M-1$ and
$$P':=\partial_1^{\beta_1}\cdots\partial_N^{\beta_N}P=\sum_{i=1}^Nc'_iX_i+c'_0\in{\cal S}\setminus\{0\},\quad c_0',c_1',\ldots,c_N'\in{\cal  R}^p,$$
where
$$\partial_1^{\beta_1}\cdots\partial_N^{\beta_N}P_1=\sum_{i=1}^Nc'_iX_i,\quad
\partial_1^{\beta_1}\cdots\partial_N^{\beta_N}P_2=c'_0.$$

If (case 1) the polynomials $c'_1,\ldots,c'_N$ are all in $U$, then we are done. Otherwise,
(case 2), there exists $i$ such that $c'_i$ is non-constant (its degree in $X_j$ is then $\geq p$
for some $j$). Now, $c'_i=\partial_iP'$ belongs to $({\cal R}^p\cap{\cal S})\setminus\{0\}$ and there
exists $s>0$ with $c'_i=P''{}^{p^s}$ with $P''\in{\cal S}$ which is not a $p$-th power.
We have constructed an element $P''$ of ${\cal S}$ which is not a $p$-th power, whose
degrees in $X_j$ are all strictly smaller than those of $P$ for all $j$ (if the polynomial depends on $X_j$).

We can repeat this process with $P''$ at the place of $P$ and so on. Since at each stage
we get a polynomial $P''$ with partial degrees in the $X_j$ strictly smaller than those of $P$ for all $j$
(if $P''$ depends on $X_j$), we eventually terminate with a polynomial $P$ which 
has all the partial degrees $<p$ in the indeterminates on which it depends, for which the 
case 1 holds.

As for the statement on the field $W$, we remark that we have applied to $P$ an algorithm
which constructs $G$ from $P$ applying finitely many partial derivatives and $p$-th roots extractions
successively, the only operations bringing out of the field $W$ being $p$-root extractions.
Hence, the existence of the integer $M$ is guaranteed.

\CVD

We recall that $U_0$ is the subfield of
$U$ whose elements are the $s\in U$ such that $\tau s=s$. Let $V$ be a subgroup of $U^\times$ such that $V\setminus V^p\not=\emptyset$.
\begin{Proposition}\label{propo2}
Under the hypotheses of Proposition \ref{propo1}, let us assume that 
for all $D\in V\setminus \{1\}$,
the only solution $s\in U$ of $\tau s=Ds$ is zero and that $D_1,\ldots,D_N\in V\setminus V^p$.
Then, the polynomial $G\in{\cal R}$ given by this Proposition 
is of the form $G=\sum_ic_iX_i+c_0$
with $c_1,\ldots,c_N\in U_0$ and $c_0\in U$. Moreover, if $c_i,c_j\not=0$ for $1\leq i<j\leq N$,
then $D_i=D_j$. Let ${\cal I}$ be the non-empty subset of $\{1,\ldots,N\}$ whose elements $i$ are such that $c_i\not=0$, let $D_i=D$ for all $i\in{\cal I}$. Then, 
$$c_0=\frac{\tau(c_0)}{D}+\frac{1}{D}\sum_{i\in {\cal I}}c_iB_i.$$
\end{Proposition}

\noindent \emph{Proof.} Proposition \ref{propo1} gives us a polynomial $G$ with $\widetilde{G}/G\in{\cal  R}$, of the form
$\sum_ic_iX_i+B^p$ with $c_i\in U$ not all vanishing and $B\in{\cal  R}$. Let $s\und{X}^{p\und{\lambda}}$ be a monomial of maximal degree in $B^p$.
Since $\widetilde{G}=QG$ with $Q\in U^\times$, we have $\tau s=(D_1^{\lambda_1}\cdots D_N^{\lambda_N})^{-p}Qs$. Moreover, $\tau(c_i)=D_i^{-1}Qc_i$ for all $i$.
Hence, if $i$ is such that $c_i\not=0$, $r:=s/c_i$ satisfies $\tau r=D_i(D_1^{\lambda_1}\cdots D_N^{\lambda_N})^{-p}r$. Now, $D_i(D_1^{\lambda_1}\cdots D_N^{\lambda_N})^{-p}\not=1$ (because
$D_i\in V\setminus V^p$)
and $r=0$, that is $s=0$. This shows that $B\in U$. Let us suppose that 
$1\leq i,j\leq N$ are such that $i\not=j$ and $c_i,c_j\not=0$. Let us write $r=c_i/c_j$; we have $\tau r=D_j/D_ir$, from which we deduce $r\in U_0$ in case $D_j/D_i=1$ and $r=0$ otherwise. The Proposition is proved dividing $\sum_ic_iX_i+B^p$ by $c_j$ with $j\not=0$ and by 
considering the relation $\widetilde{P}=QP$, once observed that $Q=D$.\CVD

We proceed, in the next two subsections, to prove Theorems \ref{papa1} and \ref{changyu}. We will prove the first theorem applying 
Propositions \ref{propo1} and \ref{propo2} to the field $U=\bigcup_{n\geq 0}C(t^{1/p^n})$ and then
by using the criterion \cite[Theorem 1.3.2]{ABP} and we will prove the second theorem 
applying these propositions to the field $U=\bigcup_{n\geq 0}K^{\text{alg.}}(x^{1/p^n})$ and then by using Theorem \ref{theorem:denis}.

\subsubsection{Direct proof of Theorem \ref{papa1}.}

For $\beta\in K^{\text{alg.}}$ such that $|\beta|<q^{q/(q-1)}$, we will use the formal series in $K^{\text{alg.}}((t))$
$$L_{\beta}(t)=\beta+\sum_{i=1}^\infty\frac{(-1)^{i}\beta^{q^i}}{(\theta^q-t)\cdots(\theta^{q^{i}}-t)},$$
defining holomorphic functions for $|t|<q^q$ with $L_\beta(\theta)=\log_\carlitz\beta$
(\footnote{Papanikolas uses these series in \cite{Papanikolas1}. It is also possible to work with the series $\sum_{i=0}^\infty e_\carlitz((\log_\carlitz\beta)/\theta^{i+1})t^i$.}).

We denote by $W$ one of the following fields: $K^{\text{alg.}},K_\infty^{\text{alg.}},C$.
For $f=\sum_ic_it^i\in W((t))$ and 
$n\in\ZZ$
 we write $f^{(n)}:=\sum_ic_i^{q^n}t^i\in W((t))$,
 so that $f^{(-1)}=\sum_ic_i^{1/q}t^i$. We have the functional equation 
$L_{\beta}^{(-1)}(t)=\beta^{1/q}+\frac{L_{\beta}(t)}{t-\theta}$.
The function $L_\beta$ allows meromorphic continuation 
to the whole $C$, with simple poles at the points $\theta^q,\theta^{q^2},\ldots,\theta^{q^n},\ldots$ of residue 
\begin{equation}(\log_\carlitz\beta)^{q},\frac{(\log_\carlitz\beta)^{q^2}}{D_1^q},\ldots,\frac{(\log_\carlitz\beta)^{q^n}}{D_{n-1}^q}.\ldots\label{eq:poles}\end{equation}

Let $\beta_1,\ldots,\beta_m$ be algebraic numbers with $|\beta|<q^{q/(q-1)}$,
let us write $L_i=L_{\beta_i}$ for $i=1,\ldots,m$. 
Let us also consider the infinite product $\boldsymbol{\Omega}$ in (\ref{eq:omega}),
converging everywhere to an entire holomorphic function with zeros at $\theta^q,\theta^{q^2},\ldots$, and write $L_0=-\boldsymbol{\Omega}^{-1}$, which satisfies the functional equation
 $$L_0^{(-1)}(t)=\frac{L_0(t)}{t-\theta},$$
with  
 $L_0(\theta)=\widetilde{\pi}$,
 meromorphic with simple poles at the points $\theta^q,\theta^{q^2},\ldots,\theta^{q^n},\ldots$,
with residues 
 \begin{equation}
 \widetilde{\pi}^{q},\frac{\widetilde{\pi}^{q^2}}{D_1^q},\ldots,\frac{\widetilde{\pi}^{q^n}}{D_{n-1}^q},\ldots\label{eq:polespi}\end{equation}
 We now prove the following Proposition.

\begin{Proposition}\label{propo5}
If the functions $L_0,L_1,\ldots,L_m$ are algebraically dependent over $K^{\text{alg.}}(t)$, then 
$\widetilde{\pi},\log_\carlitz\beta_1,\ldots,\log_\carlitz\beta_m$ are linearly dependent over $K$.
\end{Proposition}

\noindent\emph{Proof.} The functions $L_i$ are transcendental, since they have infinitely many poles.
Without loss of generality, we may assume that $m\geq 1$ is minimal so that for all 
$0\leq n\leq m$ the functions obtained from the family 
$(L_0,L_1,\ldots,L_m)$ discarding $L_n$
are algebraically independent over $K^{\text{alg.}}(t)$.

We now apply Propositions \ref{propo1} and \ref{propo2}.
We take $U:=\bigcup_{n\geq 0}C(t^{1/p^n})$, which is perfect,
and $\tau:U\rightarrow U$ the $q$-th root map on $C$ (inverse of the Frobenius map), such that $\tau(t)=t$; this is an $\FF_q$-automorphism.
Moreover, we take $N=m+1$, $D_1=\cdots=D_N=(t-\theta)^{-1}$,
\begin{eqnarray*}
(B_1,\ldots,B_N)&=&(0,\beta_{1}^{1/q},\ldots,\beta_{m}^{1/q}),
\end{eqnarray*}
and $V=(t-\theta)^\ZZ$.

Let ${\cal T}\subset C[[t]]$ be the subring of formal series converging for all $t\in C$ with $|t|\leq 1$, let
${\cal L}$ be its fraction field. Let $f\in{\cal L}$ be non-zero. 
A variant of Weierstrass preparation theorem (see \cite[Lemma 2.9.1]{and86})
yields a unique factorisation:
\begin{equation} \label{E:LLfactorization}
  f = \lambda \biggl(\prod_{|a|_{\infty} \leq 1} (t - a)^{\mathrm{ord}_a(f)}
  \biggr)\biggl(1 + \sum_{i=1}^\infty b_i t^i \biggr),
\end{equation}
where $0 \neq \lambda \in C$, $\sup_i |b_i| < 1$, and
$|b_i| \to 0$, the product being over a finite index set. 
Taking into account (\ref{E:LLfactorization}), it is a little exercise to show that $U_0=\bigcup_{i\geq 0}\FF_q(t^{1/p^i})$ and that for $D\in V\setminus\{1\}$, the solutions in $U$ of $f^{(-1)}=Df$
are identically zero (for this last statement, use the transcendence over $U$ of $\boldsymbol{\Omega}$).

Let $P\in{\cal R}$ be an irreducible polynomial such that $P(L_0,L_1,\ldots,L_m)=0$; we clearly have $\widetilde{P}=QP$ with $Q\in U$ and
Propositions \ref{propo1} and \ref{propo2} apply to give
$c_1(t),\ldots,c_{m}(t)\in U_0$ not all zero and $c(t)\in U$ such that
\begin{equation}c(t)=(t-\theta)c^{(-1)}(t)+(t-\theta)\sum_{i=1}^{m}c_i(t)\beta_{i}^{1/q}.
\label{eq:c}\end{equation}
We get, for all $k\geq 0$:
\begin{eqnarray}
c(t)&=&-\sum_{i=1}^mc_i(t)\left(\beta_i+\sum_{h=1}^k\frac{(-1)^{h}\beta_i^{q^h}}{(\theta^q-t)(\theta^{q^2}-t)\cdots(\theta^{q^h}-t)}\right)\label{eq:iteratedc}\\
& &+\frac{c^{(k+1)}(t)}{(\theta^q-t)(\theta^{q^2}-t)\cdots(\theta^{q^{k+1}}-t)}.\nonumber
\end{eqnarray}


We endow ${\cal L}$ with a norm $\|\cdot\|$ in the following way: if $f\in{\cal L}^\times$ factorises as in
(\ref{E:LLfactorization}), then $\|f\|:=|\lambda|$. 
Let $g$ be a positive integer. Then $\|\cdot\|$ extends in a unique way to the subfield
${\cal L}_g:=\{f:f^{p^g}\in{\cal L}\}$. If $(f_i)_{i\in\NN}$ is a uniformly convergent  sequence in ${\cal L}_g$
(on a certain closed ball centered at $0$) such that $\|f_i\|\rightarrow 0$, then $f_i\rightarrow 0$ uniformly.

We observe that there exists $g\geq 0$ such that $c(t),c_1(t),\ldots,c_m(t)\in {\cal L}_g$.
Hence $c_1(t),\ldots,c_m(t)\in\FF_q(t^{1/p^g})$ and $\|c_i\|=1$ if $c_i\not=0$. This implies that
$$\left\|\sum_{i=1}^mc_i\beta_i^{1/q}\right\|\leq\max_i\{|\beta_i^{1/q}|\}<q^{1/(q-1)}.$$
By (\ref{eq:c}), $\|c\|\leq q^{q/(q-1)}$. Indeed, two cases occur. The first case when $\|c^{(-1)}\|\leq \max_i\{|\beta_i^{1/q}|\}$; here we have $\|c\|<q^{q/(q-1)}$ because $\|c^{(-1)}\|=\|c\|^{1/q}$ by (\ref{E:LLfactorization}) and $\max_i\{|\beta_i|\}<q^{q/(q-1)}$ by hypothesis. The second case occurs when the inequality $\|c^{(-1)}\|>\max\{|\beta_1^{1/q}|,\ldots,|\beta_m^{1/q}|\}$ holds. In this case,
$\max\{\|c^{(-1)}(t-\theta)\|,\|(t-\theta)\sum_{i=1}^mc_i(t)\beta_i^{1/q}\|\}=\|c^{(-1)}(t-\theta)\|$
which yields $\|c\|=q^{q/(q-1)}$ by (\ref{eq:c}).

Going back to (\ref{eq:iteratedc}) we see that the sequence of functions $$E_h(t)=\frac{c^{(h+1)}(t)}{(\theta^q-t)(\theta^{q^2}-t)\cdots(\theta^{q^{h+1}}-t)}$$ converges uniformly in every closed ball included in $\{t:|t|<q^q\}$, as the series defining the functions $L_i$ ($i=1,\ldots,m$) do. We want to compute the limit of this sequence: we have two cases.

\medskip

\noindent\emph{First case.} If $\|c\|<q^{q/(q-1)}$, then, there exists $\epsilon>0$
such that $\|c\|=q^{(q-\epsilon)/(q-1)}$. Then, for all $h\geq 0$, $\|c^{(h+1)}\|=\|c\|^{q^{h+1}}=q^{(q^{h+2}-\epsilon q^{h+1})/(q-1)}$. On the other side:
\begin{eqnarray*}
\|(\theta^q-t)(\theta^{q^2}-t)\cdots(\theta^{q^{h+1}}-t)\| &=& |\theta|^{q+\cdots+q^{h+1}}\\
&=&q^{q(q^{h+1}-1)/(q-1)}.
\end{eqnarray*}
Hence, $$\|E_h\|=q^{\frac{q^{h+2}-\epsilon q^{h+1}}{q-1}-\frac{q^{h+2}-q}{q-1}}=q^{\frac{q-\epsilon q^{h+1}}{q-1}}\rightarrow 0,$$ which implies $E_h\rightarrow 0$ (uniformly on every ball as above).

This means that $\sum_{i=1}^mc_i(t)L_i(t)+c(t)=0$.
Let $g$ be minimal such that there exists a non-trivial linear relation as above, with $c_1,\ldots,c_m\in U_0\cap{\cal L}_g$; we claim that $g=0$. Indeed, if $g>0$, $c_1,\ldots,c_m\not\in\FF_q$ and 
there exists a non-trivial relation 
$\sum_{i=1}^md_i(t)L_i(t)^{p^g}+d(t)=0$ with $d_1,\ldots,d_m\in\FF_q[t]$ not all zero, $d(t)\in C(t)$ and $\max_i\{\deg_td_i\}$ minimal, non-zero. But letting the operator $d/dt$ act on this relation
we get a non-trivial relation with strictly lower degree because $dF^p/dt=0$, leading to a contradiction.

Hence, $g=0$ and $c_1,\ldots,c_m\in\FF_q(t)$. This also implies that $c\in C$; multiplying 
by a common denominator, we get a non-trivial relation
$\sum_{i=1}^mc_i(t)L_i(t)+c(t)=0$ with $c_1,\ldots,c_m\in\FF_q[t]$ and $c\in C(t)$.
The function $c$ being algebraic, it has finitely many poles. This means that $$\sum_{i=1}^mc_i(t)L_i(t)$$
has finitely many poles but for all $i$, $L_i$ has poles at $\theta^q,\theta^{q^2},\ldots$ with residues
as in (\ref{eq:poles}), which implies that $\sum_{i=1}^mc_i(t)L_i(t)$
has poles in $\theta^q,\theta^{q^2},\ldots$.
Since the functions $c_i$ belong to $\FF_q[t]$, they vanish only at points of absolute value $1$, and
the residues of the poles are multiples of $\sum_{i=0}^mc_i(\theta)^{q^k}(\log_\carlitz\beta_i)^{q^k}$ ($k\geq 1$) by non-zero factors in $A$.
They  all must vanish: this happens if and only if $$\sum_{i=1}^mc_i(\theta)\log_\carlitz\beta_i=0,$$ where we also observe that $c_i(\theta)\in K$;
the Proposition follows in this case.

\medskip

\noindent\emph{Second case.} 
Here we know that the sequence $E_h$ converges, but not to $0$ and we must compute its limit.
Let $\nu$ be in $C$ with $|\nu|=1$.
Then, there exists $\mu\in\FF_q^{\text{alg.}}{}^\times$, unique such that $|\nu-\mu|<1$.
Hence, if $\lambda\in C$ is such that $|\lambda|=q^{q/(q-1)}$,
there exists $\mu\in\FF_q^{\text{alg.}}{}^\times$ unique with 
\begin{equation}|\lambda-\mu(-\theta)^{q/(q-1)}|<q^{q/(q-1)}.\label{numu}\end{equation}
We have:
$$c(t)=\lambda\prod_{|a|\leq 1}\left(t^{1/p^g}-a\right)^{\ord_{a}c}\left(1+\sum_{i\geq 1}b_it^{i/p^g}\right),$$ with $\lambda\in C^\times$, the product being finite and $|b_i|<1$ for all $i$ so that $\|c\|=|\lambda|$.

Let $\mu\in\FF_q^{\text{alg.}}{}^\times$ be such that (\ref{numu}) holds, and 
write:
\begin{eqnarray}
c_1(t)&=&(\lambda-\mu(-\theta)^{q/(q-1)})\prod_{|a|\leq 1}\left(t^{1/p^g}-a\right)^{\ord_{a}c}\left(1+\sum_{i\geq 1}b_it^{i/p^g}\right),\nonumber\\
c_2(t)&=&\mu(-\theta)^{q/(q-1)}\prod_{|a|\leq 1}\left(t^{1/p^g}-a\right)^{\ord_{a}c}\left(1+\sum_{i\geq 1}b_it^{i/p^g}\right),\label{c2}
\end{eqnarray}
so that $c(t)=c_1(t)+c_2(t)$, $\|c_1\|<q^{q/(q-1)}$ and $\|c_2\|=q^{q/(q-1)}$.
For all $h$, we also write:
$$E_{1,h}(t)=\frac{c_1^{(h+1)}(t)}{(\theta^q-t)(\theta^{q^2}-t)\cdots(\theta^{q^{h+1}}-t)},\quad
E_{2,h}(t)=\frac{c_2^{(h+1)}(t)}{(\theta^q-t)(\theta^{q^2}-t)\cdots(\theta^{q^{h+1}}-t)}.$$
Following the first case, we easily check that $E_{1,h}(t)\rightarrow 0$
on every closed ball of center $0$ included in $\{t:|t|<q^q\}$. It remains to compute the limit of
$E_{2,h}(t)$.

We look at the asymptotic behaviour of the images of the factors in (\ref{c2}) under the operators
$f\mapsto f^{(n)}$, $n\rightarrow\infty$.
The sequence of functions $(1+\sum_{i\geq 1}b_it^{i/p^g})^{(n)}$ converges to $1$ for $n\rightarrow \infty$ uniformly on every closed ball as above.
Let ${\cal E}$ be the finite set of the $a$'s involved in the finite product (\ref{c2}), take $a\in {\cal E}$.
If $|a|<1$, then $a^{(n)}\rightarrow 0$ and $(t^{1/p^g}-a)^{(n)}\rightarrow t^{1/p^g}$.
If $|a|=1$, there exists $\mu_a\in\FF_q^{\text{alg.}}{}^\times$ such  that $|a-\mu_a|<1$ and 
we can find $n_a>0$ integer such that $\lim_{s\rightarrow \infty}a^{(sn_a)}= \mu_a$, whence 
$\lim_{s\rightarrow \infty}(t^{1/p^g}-a)^{(sn_a)}= t^{1/p^g}-\mu_a$. 

Let us also denote by $\tilde{n}>0$
the smallest positive integer such that $\mu^{q^{\tilde{n}}}=\mu$.
Let $N$ be the lowest common multiple 
of $\tilde{n}$ and the $n_a$'s with $a$ varying in ${\cal  E}$. Then the sequence of functions:
$$\left(\left(\prod_{|a|\leq 1}\left(t^{1/p^g}-a\right)^{\ord_{a}c}\right)\left(1+\sum_{i\geq 1}b_it^{i/p^g}\right)\right)^{(Ns)},\quad s\in\NN$$
converges to a non-zero element $Z\in\FF_q^{\text{alg.}}(t^{1/p^g})$.

For $n\in\NN$, let us write:
$$V_n(t):=\mu^{q^n}\frac{(-\theta)^{q^{n+1}/(q-1)}}{(\theta^{q}-t)(\theta^{q^{2}}-t)\cdots(\theta^{q^{n+1}}-t)}.$$
We have:
\begin{eqnarray*}
(-\theta)^{q/(q-1)}\prod_{i=1}^{n+1}\left(1-\frac{t}{\theta^{q^i}}\right)^{-1}&=&(-1)^{q/(q-1)}\theta^{q/(q-1)}\theta^{(q+\cdots+q^{n+1})}\prod_{i=1}^{n+1}(\theta^{q^i}-t)^{-1}\\
&=&(-1)^{q/(q-1)}\theta^{q^{n+2}/(q-1)}\prod_{i=1}^{n+1}(\theta^{q^i}-t)^{-1}.
\end{eqnarray*}
Hence, $\lim_{n\rightarrow\infty}\theta^{q/(q-1)}/((\theta^q-t)(\theta^{q^2}-t)\cdots(\theta^{q^{n+1}}-t))^{-1}=\boldsymbol{\Omega}(t)^{-1}$ 
from which we deduce that $\lim_{s\rightarrow\infty}E_{2,sN}(t)=c_0(t)L_0(t)$ with $c_0\in\FF_q^{\text{alg.}}{}^\times(t^{1/p^g})$.
We have proved that for some $c_1,\ldots,c_m\in\FF_q(t^{1/p^g}),c_0\in\FF_q(t^{1/p^g})^\times$ and $c\in K^{\text{alg.}}(t^{1/p^g})$, $\sum_{i=0}^mc_iL_i+c=0$. Applying the same 
tool used in the first case we can further prove that in fact, $g=0$. If $c_0$ is not defined over
$\FF_q$, then applying the operator $f\mapsto f^{(-1)}$ we get another non-trivial relation
$c'_0+\sum_{i=1}^mc_iL_i=c'$ with $c'\in K^{\text{alg.}}(t)$ and $c'_0\in\FF_q^\times(t)$
not equal to $c_0$; subtracting it from the former relation yields $L_0\in K^{\text{alg.}}(t)$ which is impossible
since $\boldsymbol{\Omega}$ is transcendental over $C(t)$. Hence $c_0$ belongs to $\FF_q(t)$ too.
Multiplying by a common denominator in $\FF_q[t]$ and applying arguments of the first 
case again (by using the explicit computation of the residues of the poles of $L_0$
at $\theta^q,\theta^{q^2},\ldots$), we find a non-trivial relation $c_0(\theta)\widetilde{\pi}+\sum_{i=1}^mc_i(\theta)\log_\carlitz\beta_i=0$.\CVD

\noindent\emph{Proof of Theorem \ref{papa1}}. 
If $\ell\in C$ is such that $e_\carlitz(\ell)\in K^{\text{alg.}}$, then there exist $a,b\in A$, $\beta\in K^{\text{alg.}}$ with $|\beta|<q^{q/(q-1)}$ such that 
$\ell=a\log_\carlitz\beta+b\widetilde{\pi}$. This well known property (also used in \cite{Papanikolas1}, see Lemma 7.4.1), together 
with Theorem 3.1.1 of \cite{ABP}, implies Theorem \ref{papa1}.\CVD


\subsubsection{Direct proof of Theorem \ref{changyu}.\label{directproof}}

Let $s\geq 1$ be an integer and let $\boldsymbol{\text{Li}}_n$ denote the $s$-th Carlitz's {\em polylogarithm} by:
$$\boldsymbol{\text{Li}}_s(z)=\sum_{k=0}^\infty\frac{(-1)^{ks}z^{q^k}}{([k][k-1]\cdots[1])^s},$$
so that $\boldsymbol{\text{Li}}_1(z)=\log_{\carlitz}(z)$ 
(the series $\boldsymbol{\text{Li}}_s(z)$ converges for $|z|<q^{sq/(q-1)}$).

For $\beta\in K^{\text{alg.}}\cap K_\infty$ such that $|\beta|<q^{sq/(q-1)}$ (a discussion about this hypothesis follows in \ref{differentdeformation}), we will use as in \cite{denis} the series
$$F_{s,\beta}(x)=\widetilde{\beta}(x)+\sum_{i=1}^\infty\frac{(-1)^{is}\widetilde{\beta}(x)^{q^i}}{(x^q-\theta)^s\cdots(x^{q^{i}}-\theta)^s},$$ where $\widetilde{\beta}(x)$ is the formal series in $\FF_q((1/x))$ obtained 
from the formal series of $\beta\in\FF_q((1/\theta))$ by replacing $\theta$ with $x$, an independent indeterminate. 

Let us assume that $x\in C$, with $|x|>1$. We have, for $i$ big enough,
$$\left|\frac{\widetilde{\beta}(x)^{q^i}}{(x^q-\theta)^s\cdots(x^{q^{i}}-\theta)^s}\right|=|x|^{q^i\deg_\theta\beta-q\frac{q^i-1}{q-1}},$$ so that the series $F_{s,\beta}(x)$ converges for $|x|>1$ provided that 
$|\beta|<q^{\frac{q}{q-1}}$ and $x$ is not of the form $\theta^{1/q^i}$.

We have the functional equations: 
$$F_{s,\beta}(x^q)=(x^q-\theta)^s(F_{s,\beta}(x)-\widetilde{\beta}(x)),$$ moreover, 
$$F_{s,\beta}(\theta)=\beta_j+\sum_{i=1}^\infty\frac{(-1)^i\beta^{q^i}}{(\theta^q-\theta)^s\cdots(\theta^{q^{i}}-\theta)^s}=\boldsymbol{\text{Li}}_s(\beta).$$
Therefore, these series define holomorphic functions for $|x|>q^{1/q}$ an allow 
meromorphic continuations to the open set $\{x\in C,|x|>1\}$, with poles at the points $\theta^{1/q^i}$.
We have ``deformed" certain Carlitz's logarithms and got in this way Mahler's functions
(except that the open unit disk is replaced with the complementary of the closed unit disk, but changing 
$x$ to $x^{-1}$ allows us to work in the neightbourhood of the origin). 

Let ${\cal J}$ be a finite non-empty subset of $\{1,2,\ldots\}$ such that if $n\in{\cal  J}$, $p$ does not divide $n$. Let us consider, for all $s\in{\cal J}$, an integer $l_s\geq 1$ and elements $\beta_{s,1},\ldots,\beta_{s,l_s}\in K^{\text{alg.}}\cap K_\infty$ with $|\beta_{s,i}|<q^{qs/(q-1)}$ ($i=1,\ldots,l_s$).
We remark that if $s$ is divisible by $q-1$ then, for all $r>q^{1/q}$ the product:
$$(-x)^{sq/(q-1)}\prod_{i=1}^\infty\left(1-\frac{\theta}{x^{q^i}}\right)^{-s}$$
converges uniformly in the region $\{x\in C,|x|\geq r\}$ to a holomorphic function $F_{s,0}(x)$, which
is the $(q-1)$-th power of a formal series in $K((1/(-x)^{1/(q-1)}))$, hence in $K((1/x))$ (compare with the function of \ref{anotherexample}). Moreover, $F_{s,0}(\theta)=\widetilde{\pi}^s$.

\begin{Proposition}\label{propo3}
If the functions $(F_{s,\beta_{s,1}},\ldots,F_{s,\beta_{s,l_s}})_{s\in{\cal J}}$ are algebraically dependent
over $K^{\text{alg.}}(x)$,  there exists $s\in{\cal  J}$ and a non-trivial  relation 
$$\sum_{i=1}^{l_s}c_iF_{s,\beta_{s,i}}(x)=f(x)\in K^{\text{alg.}}(x)$$ with $c_1,\ldots,c_{l_s}\in K^{\text{alg.}}$
if $q-1$ does not divide $s$, or a non-trivial relation:
$$\sum_{i=1}^{l_s}c_iF_{s,\beta_{s,i}}(x)+\lambda F_{s,0}(x)=f(x)\in K^{\text{alg.}}(x)$$ with $c_1,\ldots,c_{l_s},\lambda\in K^{\text{alg.}}$
if $q-1$ divides $s$.
In both cases, non-trivial relations can be found with $c_1,\ldots,c_{l_s},\lambda \in A$.
\end{Proposition}

\noindent\emph{Proof.} 
Without loss of generality, we may assume that ${\cal J}$ is minimal so that for all 
$n\in{\cal  J}$ and $i\in\{1,\ldots,l_n\}$ the functions obtained from the family 
$(F_{s,\beta_{s,1}},\ldots,F_{s,\beta_{s,l_s}})_{s\in{\cal J}}$ discarding $F_{n,\beta_{n,1}}$
are algebraically independent over $K^{\text{alg.}}(x)$.

We want to apply Propositions \ref{propo1} and \ref{propo2}.
We take $U:=\bigcup_{n\geq 0}K^{\text{alg.}}(x^{1/p^n})$,
and $\tau:U\rightarrow U$ the identity map on $K^{\text{alg.}}$ extended to $U$ so that $\tau(x)=x^q$.
We also take:
\begin{eqnarray*}
(X_1,\ldots,X_N)&=&(Y_{s,1},\ldots,Y_{s,l_s})_{s\in{\cal  J}},\\
(D_1,\ldots,D_N)&=&(\underbrace{(x^q-\theta)^s,\ldots,(x^q-\theta)^s}_{l_s\mathrm{ times }})_{s\in{\cal  J}},\\
(B_1,\ldots,B_N)&=&(\widetilde{\beta_{s,1}}(x),\ldots,\widetilde{\beta_{s,l_s}(x)})_{s\in{\cal  J}}.
\end{eqnarray*}
We take $V=(x^q-\theta)^\ZZ$.
We have $U_0=K^{\text{alg.}}$ and for $D\in V\setminus\{1\}$, the solutions of $f(x^q)=Df(x)$
are identically zero as one sees easily writing down a formal power series for a solution $f\in U$.

Let $P\in{\cal R}$ be an irreducible polynomial such that $P((F_{s,\beta_{s,1}},\ldots,F_{s,\beta_{s,l_s}})_{s\in{\cal J}})=0$; we clearly have $\widetilde{P}=QP$ with $Q\in U$ and
Propositions \ref{propo1} and \ref{propo2} apply. They give $s\in{\cal J}$,
$c_1,\ldots,c_{l_s}\in K^{\text{alg.}}$ not all zero and $c_0\in U$ such that
\begin{equation}c_0(x)=\frac{c_0(x^q)}{(x^q-\theta)^s}-\frac{1}{(x^q-\theta)^s}\sum_{i=1}^{l_s}c_i\widetilde{\beta_{s,i}}(x).
\label{eq:c0}\end{equation}
We now inspect this relation in more detail. To ease the notations, we write $l_s=m$ and
$F_i(x):=F_{s,\beta_{s,i}}$ ($i=1,\ldots,m$).
Since $\widetilde{\beta}(x)^q=\widetilde{\beta}(x^q)$ for all $\beta\in K$, from (\ref{eq:c0}) we get, for all $k\geq 0$:
\begin{eqnarray}
c_0(x)&=&-\sum_{i=1}^mc_i\left(\widetilde{\beta_i}(x)+\sum_{h=1}^k\frac{(-1)^{hs}\widetilde{\beta_i}(x)^{q^h}}{((x^q-\theta)(x^{q^2}-\theta)\cdots(x^{q^{h+1}}-\theta))^s}\right)\label{eq:iterated}\\
& &+\frac{c_0(x^{q^{k+1}})}{((x^q-\theta)(x^{q^2}-\theta)\cdots(x^{q^{k+1}}-\theta))^s}.\nonumber
\end{eqnarray}
By Proposition \ref{propo1}, there exists $M>0$ such that $c_0(x)^{q^M}\in K^{\text{alg.}}(x)$, which implies that $c_0(x^{q^M})\in K^{\text{alg.}}(x)$.
By equation (\ref{eq:iterated}) we see that $c_0(x)\in K^{\text{alg.}}(x)$.
 
We write $c_0(x)=\sum_{i\geq i_0}d_ix^{-i}$ with $d_i\in K^{\text{alg.}}$. The sequence $(|d_i|)_i$
is bounded; let $\kappa$ be an upper bound. If $x\in C$ is such that $|x|\geq r>q^{1/q}$ with $r$ independent on $x$, then
$|c_0(x)|=\sup_{i}|d_i||x|^{-i}\leq \kappa\sup_i|x|^{-i}\leq \kappa|x|^{\deg_xc_0}$.
Moreover, for $|x|>r$ with $r$ as above, $|x|^{q^s}>|\theta|=q$ for all $s\geq 1$ so that
 $|x^{q^s}-\theta|=\max\{|x|^{q^s},|\theta|\}=|x|^{q^s}$. Hence we get:
 $$|(x^q-\theta)(x^{q^2}-\theta)\cdots(x^{q^{k+1}}-\theta)|=|x|^{q+q^2+\cdots+q^{h+1}}=|x|^{\frac{q(q^{k+1}-1)}{q-1}}.$$
 
 Let us write:
$$R_k(x):=\frac{c_0(x^{q^{k+1}})}{((x^q-\theta)(x^{q^2}-\theta)\cdots(x^{q^{k+1}}-\theta))^s}.$$
We have, for $|x|\geq r>q^{1/q}$ and for all $k$:
\begin{equation}|R_k(x)|\leq \kappa|x|^{q^{k+1}\deg_xc_0-\frac{sq(q^{k+1}-1)}{q-1}}.\end{equation}

Since $|\beta_i|<q^{sq/(q-1)}$, $\deg_x\widetilde{\beta_i}(x)<sq/(q-1)$ for all $i$.
 In (\ref{eq:c0}) we have two cases: one if $\deg_x(c_0(x^q)/(\theta-x^q)^s)\leq\max_i\{\deg_x\widetilde{\beta_i}\}$, one if $\deg_x(c_0(x^q)/(\theta-x^q)^s)>\max_i\{\deg_x\widetilde{\beta_i}\}$.
 In the first case we easily see that $\deg_xc_0<sq/(q-1)$ (notice that $\deg_xc_0(x^q)=q\deg_xc_0$). In the second case, $\deg_x c_0=q\deg_xc_0-sq $ which implies $\deg_xc_0=sq/(q-1)$.

\medskip

\noindent\emph{First case.}
Here, there exists $\epsilon>0$ such that $\deg_xc_0=(sq-\epsilon)/(q-1)$. 
We easily check (assuming that $|x|\geq r>q^{1/q}$):
\begin{eqnarray*}
|R_k(x)|&\leq & \kappa|x|^{\frac{sq-\epsilon}{q-1}q^{k+1}-\frac{sq(q^{k+1}-1)}{q-1}}\\
&\leq&\kappa|x|^{\frac{sq-\epsilon q^{k+1}}{q-1}}
\end{eqnarray*}
and the sequence of functions $(R_k(x))_k$ converges uniformly to zero in the domain
$\{x,|x|\geq r\}$ for all $r>q^{1/q}$. Letting $k$ tend to infinity in (\ref{eq:iterated}), we find
$\sum_ic_iF_i(x)+c_0(x)=0$; that is what we expected.

\medskip

\noindent\emph{Second case.} In this case, the sequence $|R_k(x)|$ is bounded 
but does not tend to $0$. Notice that this case does not occur if $q-1$ does not divide $s$, because
$c_0\in K^{\text{alg.}}(x)$ and its degree is a rational integer. Hence we suppose that $q-1$ divides $s$. 

Let us write:
$$c_0(x)=\lambda x^{sq/(q-1)}+\sum_{i>sq/(q-1)}d_ix^i,$$ with $\lambda\in K^{\text{alg.}}{}^\times$. We have $$\lim_{k\rightarrow\infty}\frac{\sum_{i>sq/(q-1)}d_ix^{q^ki}}{((x^q-\theta)(x^{q^2}-\theta)\cdots(x^{q^{k+1}}-\theta))^s}=0$$ (uniformly on $|x|>r>q^{1/q}$),
as one verifies following the first case. 

For all $k\geq 0$ we have:
\begin{eqnarray*}
(-x)^{sq/(q-1)}\prod_{i=1}^{k+1}\left(1-\frac{\theta}{x^{q^i}}\right)^{-s}&=&(-1)^{sq/(q-1)}x^{sq/(q-1)}x^{s(q+\cdots+q^{k+1})}\prod_{i=1}^{k+1}(x^{q^i}-\theta)^{-s}\\
&=&(-1)^{sq/(q-1)}x^{sq^{k+2}/(q-1)}\prod_{i=1}^{k+1}(x^{q^i}-\theta)^{-s}.
\end{eqnarray*}
Hence we have $\lim_{k\rightarrow\infty}\lambda x^{sq/(q-1)}/((x^q-\theta)(x^{q^2}-\theta)\cdots(x^{q^{k+1}}-\theta))^s=\lambda F_{s,0}(x)$ and 
$\sum_ic_iF_i(x)+\lambda F_{s,0}(x)+c_0(x)=0$.

We now prove the last statement of the Proposition: this follows from 
an idea of Denis. The proof is the same in both cases and we work with the first only.
There exists $a\geq 0$ minimal such that the $p^a$-th powers of $c_1,\ldots,c_{l_s}$
are defined over the separable closure $K^{\text{sep}}$ of $K$. The trace 
$K^{\text{sep}}\rightarrow K$ can be extended to formal series $K^{\text{sep}}((1/x))\rightarrow K((1/x))$; its image does not vanish. 
We easily get, multiplying by a denominator in $A$,
a non-trivial relation
$$\sum_ib_iF_i(x)^{q^a}+b_0(x)=0$$ with $b_i\in A$ and $b_0(x)\in K^{\text{alg.}}(x)$.
If the coefficients $b_i$ are all in $\FF_q$, this relation is the $p^a$-th power
of a linear relation as we are looking for. If every relation has at least one of the coefficients $b_i$ not
belonging to $\FF_q$, the one with $\max_i\{\deg_\theta b_i\}$ and $a$ minimal
has in fact $a=0$
(otherwise, we apply the operator $d/d\theta$ to find one with smaller degree,
because $dg^{p^a}/d\theta=0$ if $a>0$).\CVD

The following proposition reproduces Denis' criterion
of algebraic independence in \cite{Denis2,Denis1}. It follows immediately from Theorem \ref{theorem:denis}.

\begin{Proposition}\label{theo_denis2} Let $L\subset K^{\text{alg.}}$ be a finite extension of $K$.
We consider $f_1,\ldots,f_m$ holomorphic functions in a domain $|x|>r\geq 1$
with Taylor's expansions in $L((1/x))$. Let us assume that there exist elements
$a_i,b_i\in L(x)$ ($i=1,\ldots,m$) such that $$f_i(x)=a_i(x)f_i(x^q)+b_i(x),\quad 1\leq i\leq m.$$
Let $\alpha$ be in $L$, $|\alpha|>r$, such that for all $n$, $\alpha^{q^n}$ 
is not a zero nor a pole of any of the functions $a_i,b_j$.

If the series $f_1,\ldots,f_m$ are algebraically independent over $K^{\text{alg.}}(x)$, then the values $$f_1(\alpha),\ldots,f_m(\alpha)$$ are algebraically independent over $K$.
\end{Proposition}

The next step is the following Proposition.

\begin{Proposition}\label{propo4}
If the numbers $(\boldsymbol{\text{Li}}_{s}(\beta_{s,1}),\ldots,\boldsymbol{\text{Li}}_{s}(\beta_{s,l_s}))_{s\in{\cal J}}$ are algebraically dependent
over $K^{\text{alg.}}$,  there exists $s\in{\cal  J}$ and a non-trivial linear relation 
$$\sum_{i=1}^{l_s}c_i\boldsymbol{\text{Li}}_{s}(\beta_{s,i})=0$$ with $c_1,\ldots,c_{l_s}\in A$.
If $q-1$ does not divide $s$, or a non-trivial relation:
$$\sum_{i=1}^{l_s}c_i\boldsymbol{\text{Li}}_{s}(\beta_{s,i})+\lambda\widetilde{\pi}^s=0$$ with $c_1,\ldots,c_{l_s},\lambda\in A$
if $q-1$ divides $s$.
\end{Proposition}
\noindent\emph{Proof.} 
By Proposition \ref{theo_denis2}, the functions $F_{s,i}$ ($s\in{\cal J},1\leq i\leq l_s$) are algebraically 
dependent over $K^{\text{alg.}}(t)$. Proposition \ref{propo3} applies and 
gives $s\in{\cal J}$ as well as a non-trivial linear dependence relation.
If $q-1$ does not divide $s$, by Proposition \ref{propo3}
there exists a non-trivial relation 
$$\sum_{i=1}^{l_s}c_iF_{s,\beta_{s,i}}(x)=f(x)\in K^{\text{alg.}}(x)$$ with $c_1,\ldots,c_{l_s}\in A$.
We substitute $x=\theta$ in this
relation:
$$\sum_{i=1}^{l_s}c_i\boldsymbol{\text{Li}}_s(\beta_{s,i})=f\in K^{\text{alg.}}.$$
After \cite{andth90} pp.172-176, for all $x\in C$ such that $|x|<q^{qs/(q-1)}$,
there exist $$v_1(x),\ldots,v_{s-1}(x)\in C$$ such that $$\left(\begin{array}{c}
0\\ \vdots \\ 0 \\ x
\end{array}\right)=\exp_s\left(\begin{array}{c}
v_1(x)\\ \vdots \\ v_{s-1}(x) \\ \boldsymbol{\text{Li}}_s(x)
\end{array}\right),$$ $\exp_s$ being the exponential function associated to the 
$s$-th twist of Carlitz's module. Moreover:
$$\exp_s\left(\begin{array}{c}
c_j v_1(\beta_{s,j})\\ \vdots \\ c_j v_{s-1}(\beta_{s,j}) \\ c_j\boldsymbol{\text{Li}}_s(\beta_{s,j})
\end{array}\right)=\phi_\carlitz^{\otimes s}(c_j)\left(\begin{array}{c}
0\\ \vdots \\ 0 \\ \beta_{s,j}
\end{array}\right)\in (K^{\text{alg.}})^s,\quad j=1,\ldots,n_s,$$ where
$\phi_\carlitz^{\otimes s}(c_j)$ denotes the action of
the 
$s$-th tensor power of Carlitz's module. 
By $\FF_q$-linearity, there exist numbers $w_1,\ldots,w_{s-1}\in C$ such that
$$\exp_s\left(\begin{array}{c}
w_1\\ \vdots \\ w_{s-1} \\ c
\end{array}\right)\in (K^{\text{alg.}})^s.$$ 
Yu's sub-$t$-module Theorem (in \cite{yu97}) implies the following 
analogue of Hermite-Lindemann's Theorem.
Let $G=(\GG_a^s,\phi)$
be a regular $t$-module with exponential function
$e_\phi$,  with $\phi(g)=a_0(g)\tau^0+\cdots$, for all $g\in A$.
Let $u\in C^s$ be such that $e_\phi(u)\in\GG_a^s(K^{\text{alg.}})$.
Let $V$ the smallest vector subspace of $C^s$ containing $u$, defined over $K^{\text{alg.}}$,
stable by multiplication by $a_0(g)$ for all $g\in A$. 
Then the $\FF_q$-subspace $e_\phi(V)$
of $C^s$ equals $H(C)$ with $H$ sub-$t$-module of $G$.

This result with $G$ the $s$-th twist of Carlitz's module and $e_\phi=\exp_s$
 implies the vanishing of $c$ and the $K$-linear dependence of the numbers 
$$\boldsymbol{\text{Li}}_s(\beta_{s,1}),\ldots,\boldsymbol{\text{Li}}_s(\beta_{s,l_s}).$$

If $q-1$ divides $s$ then by Proposition \ref{propo3}
there exists a non-trivial relation 
$$\sum_{i=1}^{l_s}c_iF_{s,\beta_{s,i}}(x)+\lambda F_{s,0}(x)=f(x)\in K^{\text{alg.}}(x)$$ with $c_1,\ldots,c_{l_s},\lambda\in K$. We substitute $x=\theta$ in this
relation:
$$\sum_{i=1}^{l_s}c_i\boldsymbol{\text{Li}}_s(\beta_{s,i})+\lambda\widetilde{\pi}^s=f\in K^{\text{alg.}}.$$
The Proposition follows easily remarking that,
after \cite{andth90} again, there exist $$v_1,\ldots,v_{s-1}\in C$$ such that $$\left(\begin{array}{c}
0\\ \vdots \\ 0 \\ 0
\end{array}\right)=\exp_s\left(\begin{array}{c}
v_1\\ \vdots \\ v_{s-1} \\ \widetilde{\pi}^s
\end{array}\right).$$\CVD

\noindent\emph{Proof of Theorem \ref{changyu}.}
To deduce Theorem \ref{changyu} from Proposition \ref{propo4} we quote Theorem 3.8.3 p. 187
of Anderson-Thakur in \cite{andth90} and proceed as in \cite{ChangYu}.
For all $i\leq nq/(q-1)$ there exists $h_{n,i}\in A$ such that if we set
$$P_n:=\sum_i\phi_{\carlitz}^{\otimes n}(h_{n,i})\left(\begin{array}{c}0\\ \vdots \\ 0 \\ \theta^i\end{array}
\right),$$ then the last coordinate $P_n$ is equal to $\Gamma(n)\zeta(n)$ (where 
$\Gamma(n)$ denotes Carlitz's arithmetic Gamma function).
Moreover, there exists $a\in A\setminus\{0\}$ with $\phi_{\carlitz}^{\otimes n}(a)P_n=0$ 
if and only if $q-1$ divides $n$. This implies that
$$\Gamma(n)\zeta(n)=\sum_{i=0}^{[nq/(q-1)]}h_{n,i}\boldsymbol{\text{Li}}_n(\theta^i).$$
The numbers $h_{n,i}$ are explicitly determined in \cite{andth90}.
In particular, one has
$$\zeta(s)=\boldsymbol{\text{Li}}_s(1),\quad s=1,\ldots,q-1.$$

We apply Proposition \ref{theo_denis2} and 
Proposition \ref{propo4} 
with ${\cal J}={\cal J}^\sharp\cup\{q-1\}$, ${\cal J}^\sharp$ being the set of all the integers $n\geq 1$ with $p,q-1$ not dividing $n$, $l_{q-1}=1$, $\beta_{q-1,1}=1$, and for $s\in {\cal J}^\sharp$,
$$(\beta_1,\ldots,\beta_{l_s})=(\theta^{i_0},\ldots,\theta^{i_{m_s}}),$$
where the exponents $0\leq i_0<\cdots<i_{m_s}\leq sq/(q-1)$ are chosen so that 
$$\zeta(s)\in K\boldsymbol{\text{Li}}_s(1)+\cdots+K\boldsymbol{\text{Li}}_s(\theta^{[sq/(q-1)]})= K\boldsymbol{\text{Li}}_s(\theta^{i_0})\oplus\cdots\oplus K\boldsymbol{\text{Li}}_s(\theta^{i_{m_s}}).$$\CVD

\begin{Remarque}\label{remcarlitz}
{\em With $\beta\in K$ as above, we can identify, replacing $\theta$ with $t^{-1}$, the formal series 
$F_\beta\in C((x))$ with a formal series $F_\beta^*\in\FF_q[t][[x]]\subset\mathcal{K}[[x]]$ (as in \ref{foursteps}), over which the operator $\tau$ defined there acts. Carlitz's module $\phi_\carlitz:\theta\mapsto\theta\tau^0+\tau^1$ acts on $\FF_q[t][[x]]$ and it is easy to compute the image of $F^*_\beta$ under this action.
from this we get:
$$F_{\Phi_{\text{\carlitz}}(\theta)\beta}(x)=\theta F_\beta(x)+(x-\theta)\beta(x),$$
which implies that, for all $a\in A$, $F_{\Phi_{\text{\carlitz}}(a)\beta}(x)\in a F_\beta(x)+\FF_q(\theta,x)$. In some sense, the functions $F_\beta$ are ``eigenfunctions" of the Carlitz module (a similar property holds 
for the functions $L_\beta$ and the functions (\ref{eq:slambda}), which also have Mahler's functions as counterparts.}
\end{Remarque}

\subsubsection{Final remarks.\label{differentdeformation}}

The fact that we could obtain Theorems \ref{papa1} and \ref{changyu} in a direct way 
should not induce a false optimism about Mahler's approach to algebraic independence; the matrices
of the linear $\tau$-difference equation systems involved are diagonal and we benefitted of this very special situation. In the general case, 
it seems more difficult to compute the transcendence degree of the field generated 
by solutions $\underline{f}={}^t(f_1,\ldots,f_m)\in\mathcal{K}((x))$ of a system like:
$$\underline{f}(x^d)=\mathcal{A}(x)\cdot\underline{f}(x)+\underline{b}(x)$$ (see for example, the difficulties
encountered in
\cite[Section 5.2]{Ni}).
One of the reasons is that tannakian approach to this kind of equation is, the time being, not yet explored.

This point of view should be considered since it has been very successful in the context of
$t$-motives as in Papanikolas work \cite{Papanikolas1}, which is fully compatible with Galois' approach. We could expect, once the tannakian 
theory of Mahler's functions is developed enough, to reach more general results by computing dimensions of motivic Galois groups (noticing the advantageous fact that the field of constants is here algebraically closed).

However, there is an important question we shall deal with: {\em is any ``period" of 
a trivially analytic $t$-motive (in the sense of Papanikolas in \cite{Papanikolas1}) a Mahler's value?}
For example, in \ref{directproof}, we made strong restrictions on the $\beta\in K^{\text{alg.}}$
so that finally, Denis Theorem in \cite{denis} is weaker than Papanikolas Theorem \ref{papa1}.
{\em Is it possible to avoid these restrictions in some way?}

We presently do not have a completely satisfactory answer to this question, but there seem to be some elements in favour of a positive answer. We will explain this in the next few lines.

The method in \ref{directproof} of deforming Carlitz logarithms $\log_\carlitz(\beta)$ into Mahler's functions 
requires that the $\beta$'s in $K^{\text{alg.}}$ correspond to analytic functions
at infinity. If $\beta=\sum_ic_i\theta^{-i}$ lies in $K_\infty=\FF_q((1/\theta))$, the series $\widetilde{\beta}(x):=\sum_ic_ix^{-i}$
converges for $x\in C$ such that $|x|>1$ and $\beta=\widetilde{\beta}(\theta)$.

This construction still works in the perfect closure of
the maximal tamely ramified extension $F$ of $\FF_q((x^{-1}))$ but cannot be followed easily for general $\beta\in K^{\text{alg.}}\setminus K_\infty$. Artin-Schreier's polynomial $X^p-X-\theta$
does not split over $F$. Hence, if $\xi\in K^{\text{alg.}}$ is a root of this polynomial (it has absolute value $|\xi|=q^{1/p}$), the construction fails with the presence of divergent series.

Let us consider $\ell_1=\log_\carlitz(\xi)$ and $\ell_2=\log_\carlitz(\theta)=\log_\carlitz(\xi^p)-\ell_1$. It is easy to show that 
$\ell_1,\ell_2$ are $K$-linearly independent. By Theorem \ref{papa1}, $\ell_1,\ell_2$ are algebraically independent
over $K$. The discussion above shows that it is virtually impossible to apply Proposition \ref{theo_denis2}
with the base point $\alpha=\theta$. 

We now show that it is possible to modify the arguments of \ref{directproof} and apply Proposition 
\ref{theo_denis2} with the base point $\alpha=\xi$.
 
To do so, let us consider, for $\beta\in K$, the formal series
$$\widetilde{F}_\beta(x)=\widetilde{\beta}(x)+\sum_{n=1}^\infty(-1)^n\frac{\widetilde{\beta}(x^{q^n})}{\prod_{j=1}^n(x^{q^{j+1}}-x^{q^j}-\theta)}.$$
If $|\beta|<q^{\frac{q^2}{q-1}}$ and if $|x|>1$, $x\not\in\{\xi^{1/q^j}+\lambda,j\geq 1,\lambda\in\FF_q\}$, these series 
converge. In particular, under the condition on $|\beta|$ above,
they all converge at $x=\xi$ since they define holomorphic functions on the domain
$\{x\in C,|x|>q^{1/pq}\}$, which contains $\xi$.
More precisely, the value at $x=\xi$ is:
$$\widetilde{F}_\beta(\xi)=\log_{\carlitz}(\beta(\xi)).$$

We have the functional equations $$\widetilde{F}_\beta(x^q)=(\theta-x^{q^2}+x^q)(\widetilde{F}_\beta(x)-\beta(x))$$
which tells us that the functions $\widetilde{F}_\beta$ 
define meromorphic functions in the open set $\{x\in C,|x|>1\}$. With all these observations,
it is a simple exercise to apply Proposition \ref{theo_denis2} with $\alpha=\xi$ and show 
the algebraic independence of $\ell_1,\ell_2$.

The reader can extend these computations and show the algebraic independence 
of other logarithms of elements of $K^{\text{alg.}}$. However, the choice of the base point $\alpha$
has to be made cleverly, and there is no general recipe yet. Here, the occurrence of
Artin-Schreier extensions is particularly meaningful since it is commonly observed that
every finite normal extension of $\FF_q((1/\theta))$ is contained in a finite tower of Artin-Schreier 
extensions of $\FF_q((1/\theta^{1/n}))$ for some $n$ \cite[Lemma 3]{Kedlaya}.

\begin{Remarque}{\em 
Just as in Remark \ref{remcarlitz},
the action of Carlitz's module yields the following formula
$$\widetilde{F}_{(\theta^q-\theta)\beta+\beta^q}(x)=\theta\widetilde{F}_\beta(x)+(x^q-x-\theta)\widetilde{\beta}(x).$$
Since there are natural ring isomorphisms
$K[X]/(X^p-X-\theta)\cong\FF_q[\xi]\cong A$, it could be interesting to see 
if there is some $\FF_q$-algebra homomorphism $A\rightarrow C[[\tau]]$ of which 
the functions $\widetilde{F}_\beta$ are ``eigenfunctions".}\end{Remarque}

\end{document}